\documentclass[final,oneeqnum,leqno,11pt]{siamltexmm}
\usepackage{amsmath,amsfonts,amsbsy,amssymb,graphicx,psfrag,epsfig,multicol,color}
\usepackage{showkeys}
\def\blue{\textcolor{blue}}
\def\green{\textcolor[rgb]{0,0.5,0}}
\def\red{\textcolor[rgb]{0.75,0,0}}
\def\mag{\textcolor[rgb]{0.8,0,0.6}}

\def\gray{\textcolor[rgb]{0.1,0.1,0.1}}

\newtheorem{remark}{\it Remark}
\newtheorem{hyp}{H\!}

\def\bqt{\begin{quote}} \def\eqt{\end{quote}}
\def\beq{\begin{equation}}  \def\eeq{\end{equation}}
\def\beqn{\begin{eqnarray}} \def\eeqn{\end{eqnarray}}
\def\beqnn{\begin{eqnarray*}}   \def\eeqnn{\end{eqnarray*}}
\def\barr{\begin{array}}    \def\earr{\end{array}}
\def\bit{\begin{itemize}}   \def\eit{\end{itemize}}
\def\ben{\begin{enumerate}} \def\een{\end{enumerate}}
\def\bc{\begin{center}}     \def\ec{\end{center}}
\def\btab{\begin{tabular}}  \def\etab{\end{tabular}}
\def\btabe{\begin{table}}  \def\etabe{\end{table}}
\def\bpic{\begin{picture}}  \def\epic{\end{picture}}
\def\bfig{\begin{figure}}   \def\efig{\end{figure}}
\def\qu{\quad}
\def\ms{\medskip}
\def\FNS{\footnotesize}

\def\RR{\mathbb R}

\def\LL{\mathbb{L}}
\def\deq{\stackrel{\rm def}{=}}

\def\BT{\begin{theorem}}      \def\ET{\end{theorem}}
\def\BP{\begin{proposition}}  \def\EP{\end{proposition}}
\def\BD{\begin{definition}}   \def\ED{\end{definition}}
\def\BL{\begin{lemma}}        \def\EL{\end{lemma}}
\def\BC{\begin{corollary}}    \def\EC{\end{corollary}}
\def\BR{\begin{remark}}       \def\ER{\end{remark}}
\def\BE{\begin{example}}      \def\EE{\end{example}}
\def\BH{\begin{hyp}}          \def\EH{\end{hyp}}
\def\endproof{\hfill $\Box$\newline}
\def\proof{\par\paragraph{Proof} \ignorespaces}

\def\cor{Corollary~\ref}
\def\defi{Definition~\ref}
\def\rem{Remark~\ref}
\def\thm{Theorem~\ref}
\def\prop{Proposition~\ref}
\def\lem{Lemma~\ref}

\def\noi{\noindent}
\def\lab{\label}
\def\nn{\nonumber}
\def\eq{\eqref}
\def\lp{\left}
\def\rp{\right}

\def\bm#1{\mbox{\boldmath $#1$}}
\def\disp#1{{\displaystyle #1}}
\def\mb{\mbox}

\def\A{\mathcal{A}}\def\C{{\mathcal C}}
\def\F{\mathcal{F}}\def\H{{\mathcal H}}
\def\O{\mathcal{O}}\def\P{{\mathcal P}}\def\Q{\mathcal{Q}}
\def\R{\mathcal{R}}
\def\U{\mathcal{U}}\def\V{\mathcal{V}}\def\Z{\mathcal{Z}}

\def\Bm{\mathrm{B}}\def\Dm{\mathrm{D}}
\def\Km{\mathrm{K}}\def\Nm{\mathrm{N}}\def\Qm{\mathrm{Q}}
\def\Tm{\mathrm{T}}\def\Um{\mathrm{U}}

\def\Id{I}
\def\to{\rightarrow}
\def\rank{\mathrm{rank}}
\def\lev{\mathrm{lev}\,}
\def\dim{\mathrm{dim}\,}
\def\ker{\mathrm{ker}} \def\ran{\mathrm{range}}

\def\<{\langle} \def\>{\rangle}
\def\void{\varnothing}
\def\leq{\leqslant}
\def\geq{\geqslant}

\def\all{\forall\;}
\def\exist{\exists\;}
\def\sm{\,\mathbf{\setminus}\,}
\def\z{\{0\}}

\def\t{\tilde}
\def\wt{\widetilde}
\def\h{\hat}

\def\Fd{\mathcal{F}_d}
\def\m{\mathsf{M}} \def\n{\mathsf{N}}  \def\k{\mathsf{K}}
\def\x{\times}
\def\#{\,\sharp\,}

\def\ga{\gamma}
\def\be{\beta}
\def\vs{\varpi}
\def\vt{\omega}

\def\la{\lambda}
\def\eps{\varepsilon}
\def\ka{\kappa}
\def\si{\sigma}
\def\Si{\Sigma}
\def\th{\theta}

\def\Om{\Omega}
\def\Omb{\overline{\Omega}}
\def\Ombk{\overline{\Omega}_\k}
\def\Omm{\Om_{\max}}

\def\hu{\hat u}
\def\tu{\widetilde u}
\def\bu{\bar{u}}
\def\bsi{\bar{\sigma}}
\def\hsi{{\hat\sigma}}

\def\Ahs{A_{\hsi}} \def\Phs{\Pi_{\hsi}}
 \def\Pbs{\Pi_{\bsi}}
\def\Avs{A_{\vs}}  \def\Pvs{\Pi_{\vs}}
\def\Avt{A_{\vt}}  \def\Pvt{\Pi_{\vt}}
\def\hvt{\hu_{\vt}}
\def\hvs{\hu_{\vs}}

\def\uhs{\hu_{\hsi}}

\def\bv{\bar v}

\def\II{\mathbb{I}}
\def\IM{\II_{\m}}
\def\Im{\II_{\m-1}}
\def\Ir{\II_r}
\def\IN{\II_{\n}}
\def\Ik{\II_{\k}}
\def\RM{\RR^{\m}}
\def\RN{\RR^{\n}}
\def\RMN{\RR^{\m\x\n}}

\def\cpa{{\rm{($\,\P_{\vt}\,$)}~}}
\def\cpar{{\rm{($\,\Q_{\vt}\,$)}~}}
\def\cphsi{{\rm{($\,\P_{\hsi}\,$)}~}}

\title{Description of the minimizers of least squares regularized~ with~ $\bm{\ell_0}$-norm.\\
Uniqueness of the global minimizer}

\author{Mila Nikolova\thanks{CMLA, ENS Cachan, CNRS,  61 Av. Pr\'es. Wilson, 94230 Cachan, France
(\email{nikolova@cmla.ens-cachan.fr}). }}

\begin{document}
\maketitle
\newcommand{\slugmaster}{%
\slugger{siims}{xxxx}{xx}{x}{x--x}}

\vspace*{-4cm}

\renewcommand{\thefootnote}{\alph{footnote}}
\hspace*{5cm}\gray{\sc (in press\footnote {Received by the editors November 10, 2011; accepted for publication (in revised form) January 14, 2013;
published electronically DATE.
})}

\vspace*{4cm}

\renewcommand{\thefootnote}{\arabic{footnote}}
\setcounter{footnote}{0}
\begin{abstract}
We have an $\m\x\n$ real-valued arbitrary matrix $A$ (e.g. a dictionary) with $\m<\n$ and data $d$
describing the sought-after object with the help of $A$.
This work provides an in-depth analysis of the (local and global) minimizers of an objective function $\Fd$
combining a quadratic data-fidelity term and an $\ell_0$ penalty applied to each
entry of the sought-after solution, weighted by a regularization parameter $\be>0$.
For several decades, this objective has attracted a ceaseless effort to conceive
algorithms approaching a good minimizer.
Our theoretical contributions, summarized below, shed new light on the
existing algorithms  and can help the conception of innovative numerical schemes.
To solve the normal equation associated with any $\m$-row submatrix of $A$
is equivalent to compute a local minimizer $\hu$ of $\Fd$.
(Local) minimizers $\hu$ of $\Fd$ are strict if and only if
the submatrix, composed of those columns of $A$ whose indexes form the support of $\hu$,
has full column rank.
An outcome is that strict local minimizers of $\Fd$ are easily computed  without knowing the value of $\be$.
Each strict local minimizer is linear in data.
It is proved that $\Fd$ has global minimizers and that they are always strict.
They are studied in more details under the (standard) assumption that $\rank(A)=\m<\n$.
The global minimizers with $\m$-length support are seen to be impractical.
Given $d$, critical values $\be_\k$ for any  $\k\leq\m-1$ are exhibited such that if
$\be>\be_\k$, all global minimizers of $\Fd$ are $\k$-sparse.
An assumption on $A$ is adopted and proved to fail only on a closed negligible subset.
Then for all data $d$ beyond a closed negligible subset,
the objective $\Fd$ for $\be>\be_\k$, $\k\leq\m-1$, has a unique global minimizer
and this minimizer is $\k$-sparse.
Instructive small-size ($5\x 10$) numerical illustrations confirm the main
theoretical results.
\end{abstract}



\paragraph{\hspace*{-6mm}\textbf{\small DOI}} {\small 10.1137/11085476X}

\bc\gray{--------------------------------------------------------------------------------------------------------------------}\ec

\pagestyle{myheadings}
\thispagestyle{plain}
\markboth{MILA NIKOLOVA}{THE MINIMIZERS OF LEAST SQUARES REGULARIZED~ WITH~ $\bm{\ell_0}$-NORM}

\section{Introduction}

Let $A$ be an arbitrary matrix (e.g., a dictionary) such that
\[~A\in\RMN~~~\mb{for}~~~\m<\n~,\]
where the positive integers $\m$ and $\n$ are fixed.
Given a data vector $d\in\RM$, we consider an objective function $\Fd:\RN\to\RR$ of the form
\beqn\Fd(u)&=&\|Au-d\|^2_2+\be\|u\|_0~,~~\be>0~,\lab{fd}\\
\|u\|_0&=&\#\si(u)~, \nn\lab{rj}
\eeqn
where $u\in\RN$ contains the coefficients describing the sought-after object,
 $\be$ is a regularization parameter,
$\#$ stands  for cardinality and $\si(u)$ is the support of $u$ (i.e.,
the set of all $i\in\{1,\cdots,\n\}$ for which the $i$th entry of $u$ satisfies $u[i]\neq0$).
By an abuse of language, the penalty in~\eq{fd} is called the $\ell_0$-norm.
Define $\phi:\RR\to\RR$ by
\beq
\phi(t)\deq\lp\{\barr{lll}0&\mb{\rm if}&t=0~,\\1&\mb{\rm if}&t\neq0~.\earr\rp.\lab{phi}\eeq
Then $\disp{\|u\|_0=\sum_{i=1}^\n\phi(u[i])=\sum_{i\in\si(u)}\phi(u[i])}$, so $\Fd$ in~\eq{fd} equivalently reads
\beq
    \Fd(u)=\|Au-d\|^2_2+\be\sum_{i=1}^\n\phi(u[i]) =
    \|Au-d\|^2_2+\be\sum_{i\in\si(u)}\phi(u[i])~.\lab{fds}
\eeq
We focus on all (local and global) minimizers $\hu$ of an objective $\Fd$ of the form~\eq{fd}:
\beq\hu\in\RN~~\mb{such that}~~
\Fd(\hu)=\min_{u\in\O}\Fd(u)~,\lab{P0}\eeq
where $\O\,$ is an open neighborhood of $\hu$~.
We note that finding a global minimizer of $\Fd\,$ must be an {\em NP-hard}
computational problem \cite{DavisMallatAvellaneda,Tropp06}.

The function $\phi$ in~\eq{phi} served as a regularizer for a long time.
In the context of Markov random fields it was used by Geman and Geman in 1984 \cite{Geman84} and
Besag in 1986 \cite{Besag86}
 as a prior in MAP energies to restore labeled images. The MAP objective reads as
\beq\Fd(u)=\|Au-d\|^2_2+\be\sum_{k}\phi(D_ku)~,\lab{map}\eeq
where $D_k$ is a finite difference  operator and $\phi$ is given by~\eq{phi}.
This label-designed form is known as the Potts prior model, or as the multi-level logistic model
\cite{Besag89,Li95}.
Various stochastic and deterministic algorithms have been considered to minimize~\eq{map}.
Leclerc \cite{Leclerc89} proposed in 1989 a deterministic continuation method to restore piecewise constant images.
Robini, Lachal and Magnin \cite{Robini07} introduced the stochastic continuation approach
and successfully used it to reconstruct 3D tomographic images.
Robini and Magnin refined the method and the theory  in \cite{Robini10}.
Very recently, Robini and Reissman \cite{Robini12} gave theoretical results
relating the probability for global convergence  and the computation speed.

The problem stated in~\eq{fd} and~\eq{P0}---to (locally) minimize $\Fd$---arises
when {\em sparse\,} solutions are desired.
Typical application fields are signal and image processing,
morphologic component analysis, compression, dictionary building,
inverse problems, compressive sensing, machine learning, model selection, classification, and
subset selection, among others.
The original hard-thresholding method proposed by Donoho and Johnstone \cite{Donoho94} amounts
to\footnote{As a reminder, if $d$ are some noisy coefficients, the restored coefficients $\hu$
minimize $\|u-d\|^2+\be\|u\|_0$ and read $\hu[i]=0$ if $\lp|d[i]\rp|\leq\sqrt{\be}$
and $\hu[i]=d[i]$ otherwise. }
minimizing~$\Fd$, where  $d$ contains the coefficients of a signal or an image
expanded in a wavelet basis ($\m=\n$).
When $\m<\n$, various (usually strong) restrictions on $\|u\|_0$ (often $\|u\|_0$ is replaced by
a less irregular function)
and on $A$ (e.g., RIP-like criteria, conditions on $\|A\|$, etc.) are
needed to conceive numerical schemes approximating a minimizer of $\Fd$,
to establish local convergence and derive the asymptotic of the obtained solution.
In statistics the problem has been widely considered for subset selection,
and numerous algorithms have been designed, with limited
theoretical production, as explained in the book by Miller~\cite{Miller02}.
More recently, Haupt and Nowak \cite{Haupt06} investigate the statistical performances of the global minimizer
of~$\Fd$ and propose an iterative bound-optimization procedure.
Fan and Li \cite{Fan06} discuss a variable splitting and penalty decomposition
minimization technique for \eq{fd}, along with other approximations of the $\ell_0$-norm.
Liu and Wu \cite{Liu07} mix the $\ell_0$ and $\ell_1$ penalties, establish
some asymptotic properties of the new estimator and use
mixed integer programming aimed at global minimization.
For model selection, Lv and Fan \cite{Lv09} approximate the $\ell_0$ penalty using functions that are
concave on $\RR_+$ and prove a nonasymptotic nearly oracle property of the resultant estimator.
Thiao, Dinh, and Thi~\cite{Thiao08} reformulate the problem so that an approximate solution
can be found using difference-of-convex-functions programming.
Blumensath and Davies \cite{Davies08} propose an iterative thresholding scheme to
approximate a solution and prove convergence to a local minimizer of $\Fd$.
Lu and Zhang \cite{LuZhang10} suggest a penalty decomposition method to minimize $\Fd$.
Fornasier and Ward \cite{Fornasier10} propose an iterative thresholding algorithm for minimizing
 $\Fd$ where $\ell_0$ is replaced by a reasonable sparsity-promoting relaxation
given by $\phi(t)=\min\{|t|,1\}$; then convergence to a local minimizer is established.
In a recent paper by Chouzenoux et al.~\cite{Pesquet12}, a mixed
$\ell_2-\ell_0$ regularization is considered:
a slightly smoothed version of the objective is analyzed and
 a majorize-minimize subspace approach, satisfying a finite length property, converges to a critical point.
Since the submission of our paper,
image reconstruction methods have been designed where $\ell_0$ regularization is applied to the
coefficients of the expansion of the sought-after image in a wavelet frame~\cite{Zhang12,Dong12}:
the provided numerical results outperform $\ell_1$ regularization for a reasonable
computational cost achieved using penalty decomposition techniques.
In a general study on the convergence of descent
methods for nonconvex objectives, Attouch, Bolte, and Svaiter~\cite{Attouch11}
apply an inexact forward-backward splitting scheme to find a critical point of $\Fd$.
Several other references can be evoked, e.g., \cite{Schorr05,Gasso09}.

Even though overlooked for several decades, the objective $\Fd$ was essentially considered
from a numerical standpoint.
The motivation naturally comes from the promising applications and
the intrinsic difficulty of minimizing $\Fd$.

{\em The goal of this work is to analyze the (local and global)
minimizers $\hu$ of objectives $\Fd$ of the form~\eq{fd}.
\bit
\item We provide detailed results on the minimization problem.
\item The uniqueness of the global minimizer of $\Fd$ is examined as well.
\eit}
\noi We do not propose an algorithm.
However, our theoretical results raise salient questions about the existing algorithms
and can help the conception of innovative numerical schemes.

The minimization of $\Fd$ in~\eq{fd} might seem close to its constraint variants:
\beq \lab{ka}\barr{lll}
\mb{given } \eps\geq0,~&\mb{minimize}~~\|u\|_0&~\mb{subject to}~~\|Au-d\|^2\leq \eps~,\\
\mb{given } K\in\IM,~&\mb{minimize}~\|Au-d\|^2&~\mb{subject to}~~\|u\|_0\leq K~.
\earr\eeq
The latter problems are abundantly studied in the context of sparse recovery in different fields.
An excellent account is given in \cite{DonBruckElad09}, see also the book \cite{Mallat08}.
For recent achievements, we refer the reader to~\cite{CSR-URL}.
It is worth emphasizing that in general, {\em there is no equivalence between the problems stated in~\eq{ka}
and the minimization of $\Fd$ in~\eq{fd}} because all of these  problems are nonconvex.

\subsection{Main notation and definitions}\lab{ntd}

We recall that if $\hu$ is a {\em (local) minimizer} of $\Fd$, the value $\Fd(\hu)$ is a (local)
minimum\footnote{These two terms are often confused in the literature.}
of $\Fd$ reached at (possibly numerous) points $\hu$.
Saying that a (local) minimizer $\hu$ is {\em strict} means that
there is a neighborhood $\O\subset\RN$, containing $\hu$, such that
$\Fd(\hu)<\Fd(v)$  for any $v\in\O\sm\{\hu\}$.
So $\hu$ is an isolated minimizer.

Let $\k$ be any  positive integer. The expression
$\big\{u\in\RR^\k~:~u~~\mb{satisfying property}~~\mathfrak{P}\big\}$
designates the subset of $\RR^\k$ formed from all elements $u$ that meet $\mathfrak{P}$.
The identity operator on $\RR^\k$ is denoted by $\Id_\k$.
The entries of a vector $u\in\RR^\k$ read as $u[i]$, for any $i$.
The $i$th vector of the
canonical basis\footnote{More precisely, for any $i\in\Ik$, the vector $e_i\in\RR^\k$
is defined by $e_i[i]=1$ and $e_i[j]=0,~\all j\in\Ik\sm\{i\}$. }
 of $\RR^\k$ is denoted by $e_i\in\RR^\k$.
Given $u\in\RR^\k$ and $\rho>0$, the {\em open} ball at $u$ of radius $\rho$ with respect to the
$\ell_p$-norm for $1\leq p\leq\infty$ reads as
\[ \Bm_p(u,\rho)\deq\{v\in\RR^\k~:~\|v-u\|_p<\rho\}~.\]
To simplify the notation, the $\ell_2$-norm  is {\em systematically} denoted by
$$\|\cdot\|\deq\|\cdot\|_2~.$$
We denote by  $\Ik$ the {\em totally and strictly ordered} index set\footnote{E.g. without strict order we have
$\omega=\{1,2,3\}=\{2,1,1,3\}$ in which case the notation in \eq{as}-\eq{us} below is ambiguous.}
\beq\Ik\deq\big(\{1,\cdots,\k\},<\big)~,\lab{ik}\eeq
where the symbol $<$ stands for the natural order of the positive integers.
Accordingly, {\em any subset $\vt\subseteq\Ik$ inherits the property of being  totally and strictly ordered.}

We shall often consider the index set $\IN$.
The complement of $\vt\subseteq\IN$ in $\IN$  is denoted by
\[\vt^c=\IN\sm\vt\subseteq \IN~.\]

\BD\lab{sp}
For any $u\in\RR^\n$, the {\em support}  $\si(u)$ of $u$ is defined by
\beq\si(u)=\Big\{i\in\IN~:~ u[i]\neq 0\Big\}\subseteq \IN~.\nn\lab{s}\eeq
\ED
If $u=0$, clearly $\si(u)=\void$.

The $i$th column in a matrix $A\in\RMN$ is denoted by $a_i$.
It is {\em systematically} assumed that
\beq\mb{\framebox{$\disp{~a_i\neq0~~~~\all i\in\IN~.}$}}\lab{ao}\eeq
For a matrix $A\in\RMN$ and a vector $u\in\RN$, with any $\vt\subseteq\IN$,
we associate the {\em submatrix} $\Avt$ and the {\em subvector} $u_\vt$ given by
\beqn
\Avt&\deq&\big(a_{\vt[1]},\cdots,a_{\vt[\#\vt]}\big)\in\RR^{\m\times \#\vt}~,\lab{as}\\
u_\vt&\deq&\Big(u\big[\vt[1]\big],\cdots,u\big[\vt[\#\vt]\big]\Big)\in\RR^{\#\vt}\lab{us}~,
\eeqn
respectively,
as well as the zero padding operator $Z_\vt~:~\RR^{\#\vt}\to\RN$ that inverts~\eq{us}:
\beq u=Z_\vt\lp(u_\vt\rp)~,~~u[i]=\lp\{
\barr{ll}0&\mb{\rm if}~~i\not\in\vt~, \\
u_\vt[k]&\mb{for the unique $k$ such that}~~\vt[k]=i.\earr
\rp.\lab{zo}\eeq
Thus for $\vt=\void$ one finds
$u_{\void}=\void~~~\mb{\rm and}~~~u=Z_{\void}\lp(u_{\void}\rp)=0\in\RN~.$

Using \defi{sp} and the notation in \eq{as}-\eq{us}, for any $u\in\RN\sm\z$ we have
\beq\lab{sd}
\vt\in\IN~~\mb{\rm and}~~\vt\supseteq\si(u)~~~\Rightarrow~~~Au=\Avt u_{\vt}~.\eeq

To simplify the presentation, we adopt the following {\em definitions}
\footnote{Note that $(a)$ corresponds to the zero mapping on $\RR^0$ and
that $(b)$ is the usual definition for the rank of an empty matrix.}:
\beq\barr{lll}
(a)&& A_{\void}=[\ ]\in\RR^{\m\x 0}~,\\
(b)&&\rank\lp(A_{\void}\rp)=0~.\earr\lab{cw}\eeq
In order to avoid possible ambiguities\footnote{In the light of~\eq{as}, $A_\vt^T$ could
also mean $\lp(A^T\rp)_\vt$.}, we set
\[\Avt^T\deq\lp(\Avt\rp)^T~,\]
where the superscript $^T$ stands for transposed.
If $\Avt$ is invertible, similarly $\Avt^{-1}\deq\lp(\Avt\rp)^{-1}$.

In the course of this work, we shall frequently refer to the constrained
quadratic optimization problem stated next.

{\em Given $d\in\RM$    and $\vt\subseteq\IN$, problem \cpa reads as:
\beq~~~~~~~~~~~~~~~~~~~~~~~~~~~\framebox{\mb{$\lp\{\barr{lll}
&&\disp{\min_{u\in\RN }\|Au-d\|^2}~,\\~\\
\mb{\rm subject to}&& u[i]=0,~~\all i\in \vt^c~.\\
\earr\rp.$}}~~~~~~~~~~~~~~~~~~~~~~~~~~~~~~~~~~~~\mb{\cpa}\lab{cpa}\eeq  }
{\em Clearly, problem \cpa always admits a solution.}

The definition below will be used to evaluate the extent of some subsets and assumptions.
\BD\lab{ps}
A property (an assumption) is called  {\em generic on} $\RR^\k$ if it
holds true on a {\em dense open} subset of $\RR^\k$.
\ED

As usual,  a subset ${\mathcal{S}}\subset\RR^\k$ is said to be {\em negligible} in $\RR^\k$
if there exists $\Z\subset\RR^\k$
whose Lebesgue measure in $\RR^\k$ is $\LL^\k(\Z)=0$ and ${\mathcal{S}}\subseteq\Z$~.
If a property fails only on a negligible set, it is said to hold
 {\em almost everywhere}, meaning ``with probability one''.
\defi{ps} requires much more than {\em almost everywhere\,}. Let us explain.

{\em If a property holds true for all $v\in\RR^\k\sm{\mathcal{S}}$, where ${\mathcal{S}}\subseteq\Z\subset\RR^\k$,
$\Z$ is {\em closed in $\RR^\k$} and $\LL^\k(\Z)=0$, then this property is {\em generic} on $\RR^\k$.}
Indeed, $\RR^\k\sm\Z$ contains a {\em dense open} subset of $\RR^\k$.
So if a property is generic on $\RR^\k$, then it holds true almost everywhere on $\RR^\k$.
But the converse is false: an almost everywhere true property is not generic if the closure of its
negligible subset has a positive measure,\footnote{
There are many examples---e.g. $\Z=\{x\in[0,1]~:~x ~\mb{is rational}\}$,
then $\LL^1(\Z)=0$ and $\LL^1(\mb{closure}(\Z))=1$.}
because then $\RR^\k\sm\Z$ does not contains an open subset of $\RR^\k$.
In this sense, a generic property is stable with respect to the objects to which it applies.

{\em The elements of a set ${\mathcal{S}}\subset\RR^\k$ where a generic property fails are highly exceptional in~$\RR^\k$.}
The chance that a truly random $v\in\RR^\k$---i.e., a $v$ following a non singular
probability distribution on $\RR^\k$---comes across such an ${\mathcal{S}}$ can be ignored in practice.

\subsection{Content of the paper}
The main result in section~\ref{LM} tells us that finding a solution of \cpa for
$\vt\subset\IN$  is {\em equivalent} to computing a (local) minimizer of $\Fd$.
In section~\ref{sSM} we prove that the (local) minimizers $\hu$ of $\Fd$ are {\em strict}
if and only if the submatrix $A_{\si(\hu)}$ has full column rank.
The strict minimizers of $\Fd$ are shown to be linear in data $d$.
The importance of the $(\m-1)$-sparse strict minimizers is emphasized.
The global minimizers of $\Fd$ are studied in section~\ref{sgm}.
Their existence is proved.
They are shown to be strict for any $d$ and for any $\be>0$.
More details are provided under the standard assumption that $\rank(A)=\m<\n$.
Given $d\in\RM$, critical values $\be_\k$ for $\k\in\Im$ are exhibited such that all
global minimizers of $\Fd$ are $\k$-sparse\footnote{
As usual, a vector $u$ is said to be $\k$-sparse if $\|u\|_0\leq\k$.}
 if $\be>\be_\k$.

In section~\ref{kb}, a gentle assumption on $A$ is shown to be {\em generic} for all
$\m\x\n$ real matrices.
Under this assumption, for all data $d\in\RM$ beyond a closed negligible subset,
the objective $\Fd$ for $\be>\be_\k$, $\k\in\Im$, has a unique global minimizer and
this minimizer is $\k$-sparse.

Small size ($A$ is $5\x10$) numerical tests in section~\ref{nm}
illustrate the main theoretical results.

\section{All minimizers of $\bm{\Fd}$}\lab{LM}

\subsection{Preliminary results}\lab{pr}

First, we give some basic facts on problem \cpa as defined in \eq{cpa} that are needed for later use.
If $\vt=\void$, then  $\vt^c=\IN$, so the unique solution of \cpa is $\hu=0$.
For an arbitrary $\vt\subset\IN$  meeting $\#\vt\geq1$, \cpa amounts to minimizing
a quadratic term with respect to only $\#\vt$ components of $u$, the remaining entries being null.
This {\em quadratic} problem \cpar reads as
\beq\!~~~~~~~~~~~~~~~~~~~~~~~~~~~~~~~~~~~~~~~~~
  \disp{\min_{v\in\RR^{\#\vt}}\big\|\Avt v-d\big\|^2},~~\#\vt\geq1~,~\lab{cpr}
~~~~~~~~~~~~~~~~~~~~~~~~~~~~~~~~~~~~~~~~~~~~~~~~~\mb{\cpar}
\eeq
and it always admits a solution.
Using the zero-padding operator $Z_\vt$ in \eq{zo}, we have
\[\Big[~\hvt\in\RR^{\#\vt}~\mb{\rm solves \cpar and}~~\hu=Z_\vt\lp(\hvt\rp)~\Big]~~~\Leftrightarrow~~~
\Big[~\mb{$\hu\in\RN$ solves \cpa},~~\#\vt\geq1~\Big]~.\]
The optimality conditions for \cpar, combined with the definition in~\eq{cw}(a), give rise to
the following equivalence, which holds true for any $\vt\subseteq\IN$:
\beq \Big[~\mb{$\hu\in\RN$ solves \cpa}\Big]~~\Leftrightarrow~~
\Big[~\hvt\in\RR^{\#\vt}~~\mb{\rm solves $~\disp{\Avt^T\Avt\, v=\Avt^Td~}$ and}~~\hu=Z_\vt\lp(\hvt\rp)~\Big].
\lab{eu}\eeq
Note that $\disp{\Avt^T\Avt\, v=\Avt^Td}$ in~\eq{eu} is the normal equation associated with $\Avt\,v=d$.
The remark below shows that the optimal value of \cpa in \eq{cpa} can also be seen as an orthogonal projection problem.

\begin{remark}  \lab{ref}  \rm
Let $r\deq\rank(\Avt)$  and ${B_{\vt}}\in\RR^{\m\x r}$ be an orthonormal basis for $\ran(\Avt)$.
Then $\Avt=B_\vt H_\vt$ for a uniquely defined matrix ${H_{\vt}}\in\RR^{r\x\#\vt}$ with $\rank({H_{\vt}})=r$.
Using~\eq{eu}, we have
\beqnn \Avt^T\Avt\hvt=\Avt^Td~~\Leftrightarrow~~{H_{\vt}}^T{H_{\vt}}\hvt={H_{\vt}}^T{B_{\vt}}^Td
~~\Leftrightarrow~~{H_{\vt}}\hvt={B_{\vt}}^Td~~\Leftrightarrow~~\Avt\hvt={B_{\vt}}{B_{\vt}}^Td~.
\eeqnn
In addition, $\Pi_{\ran(\Avt)}={B_{\vt}}{B_{\vt}}^T$ is the orthogonal projector onto the subspace spanned
by the columns of $\Avt$, see e.g.~\cite{Meyer00}. The expression above combined with~\eq{eu} shows that
\[\Big[\,\mb{$\hu\in\RN$ solves }(\P_{\vt})\,\Big]~~\Leftrightarrow~~
\Big[~\hvt\in\RR^{\#\vt}~~\mb{\rm meets}~~\Avt\hvt=\Pi_{\ran(\Avt)}d~~~\mb{and}~~\hu=Z_\vs(\hvs)\,\Big].
\lab{pref}\]
\em Obviously,
$A\hu=\Avt\hvt$ is the orthogonal projection of $d$ onto the basis~$B_\vt$.\end{remark}

For $\vt\subseteq\IN$, let $\Km_\vt$ denote the vector subspace
\beq\Km_\vt\deq\big\{v\in\RN~:~v[i]=0,~\all i\in\vt^c\big\}~.\lab{km}\eeq
This notation enables problem \cpa in \eq{cpa} to be rewritten as
\beq\disp{\min_{u\in\Km_\vt}\|Au-d\|^2~.}\lab{fa}\eeq

The technical lemma below will be used in what follows.
We emphasize that its statement is {\em independent} of the vector $\hu\in\RN\sm\z$.

\BL\lab{tl} Let $d\in\RM$,  $\be>0$, and $\hu\in\RN\sm\z$ be {\em arbitrary}. For $\hsi\deq\si(\hu)$, set
\beq
\rho\deq\min\lp\{\min_{i\in\hsi}\big|\,\hu[i]\,\big|,~\frac{\be}{2\Big(\|A^T(A\hu-d)\|_1+1\Big)}\rp\}.
\lab{rho}\eeq
Then $\rho>0$.
\bit
\item[\rm(i)] For $\phi$ as defined in~\eq{phi}, we have
\beq\nn\lab{tlo}
v\in \Bm_\infty(0,\rho)~~\Rightarrow~~\sum_{i\in\IN}\phi\big(\,\hu[i]+v[i]\,\big)=
\sum_{i\in\hsi}\phi\lp(\hu[i]\rp)+\sum_{i\in\hsi^c}\phi\lp(v[i]\rp)~.
\eeq
\item[\rm(ii)] For $\Km_{\hsi}$ defined according to \eq{km}, $\Fd$ satisfies
\beq  v\in\Bm_\infty(0,\rho)\cap\big(\RN\sm\Km_{\hsi}\big)~~~\Rightarrow~~~\Fd(\hu+v)\geq\Fd(\hu) ~,
\nn\lab{pip}\eeq
where the inequality is {\em strict} whenever  $\hsi^c\neq\void$.
\eit
\EL

The proof is outlined in Appendix~\ref{dtl}.

\subsection{The (local) minimizers of $\bm{\Fd}$ solve quadratic problems}\lab{fip}~~

{\em It is worth emphasizing that no special assumptions on the matrix $A$ are adopted.}

 We begin with an easy but cautionary result.

\BL\lab{uz}
For any $d\in\RM$ and for all $\be>0$, $\Fd$ has a {\em strict} (local) minimum at $\hu=0\in\RN$.
\EL

\proof
Using the fact that $\Fd(0)=\|d\|^2\geq0$, we have
\beqn
\Fd(v)&=&\|Av-d\|^2+\be\|v\|_0=\Fd(0)+\R_d(v)~,\lab{wq}\\
\mb{where}\quad
\R_d(v)&=&\|Av\|^2-2\<v,A^Td\>+\be\|v\|_0~.\eeqn
Noticing that $\be\|v\|_0\geq\be>0$ for $v\neq0$ leads to
\[v\in \Bm_2\lp(0,\frac{\be}{2\|A^Td\|+1}\rp)\setminus\{0\}~~~\Rightarrow~~~
\R_d(v)\geq-2\|v\|\;\|A^Td\|+\be >0~.\]
Inserting this implication into~\eq{wq} proves the lemma.
\endproof

{\em For any $\be>0$ and $d\in\RM$, the sparsest strict local minimizer of $\Fd$ reads $\hu=0$.
Initialization with zero of a suboptimal algorithm should generally be a bad choice.}
Indeed, experiments have shown that such an initialization can be harmful;
see, e.g., \cite{Miller02,Davies08}.

The next proposition states a result that is often evoked in this work.

\BP\lab{cp}
Let $d\in\RM$.
Given an $\vt\subseteq\IN$, let  $\hu$ solve problem \cpa as formulated in~\eq{cpa}.
Then for any $\be>0$, the objective $\Fd$ in~\eq{fd} reaches a (local) minimum at $\hu$
and \beq\si(\hu)\subseteq\vt~,\lab{cop}\eeq
where $\si(\hu)$ is given in \defi{sp}.
\EP

\proof
Let $\hu$ solve problem \cpa, and let $\be>0$. The constraint in \cpa entails~\eq{cop}.

Consider that $\hu\neq0$, in which case for $\hsi\deq\si(\hu)$ we have $1\leq\#\hsi\leq\#\vt$.
Using the equivalent formulation of \cpa given in~\eq{km}-\eq{fa}, yields
\beq v\in\Km_\vt~~~~\Rightarrow~~~\|A(\hu+v)-d\|^2\geq\|A\hu-d\|^2~.\lab{sqi}\eeq
The inclusion in \eq{cop} is equivalent to
$\vt^c\subseteq\hsi^c~.$
Let $\Km_{\hsi}$ be defined according to~\eq{km} as well.
Then
\[\hu\in\Km_{\hsi}\subseteq\Km_\vt.\]
Combining~the latter relation with~\eq{sqi} leads to
\beq\lab{yui}
v\in\Km_{\hsi}~~~\Rightarrow~~~\|A(\hu+v)-d\|^2\geq\|A\hu-d\|^2~.
\eeq
Let $\rho$ be defined  as in~\eq{rho}  \lem{tl}.
Noticing that by \eq{phi} and \eq{km}
\beq v\in\Km_{\hsi}~~~\Rightarrow~~~\phi\lp(v[i]\rp)=0~~~\all i\in\hsi^c~, \lab{rer}\eeq
the following inequality chain is derived:
\beqnn v\in \Bm_\infty(0,\rho)\cap\Km_{\hsi}~~~\Rightarrow~~~
\Fd(\hu+v)&=&\|A(\hu+v)-d\|^2+\be\sum_{i\in\IN}\phi\lp(\hu[i]+v[i]\rp)\nn\\
\Big[\mb{\rm by \lem{tl}(i)}\Big]~~~~
    &=&\|A(\hu+v)-d\|^2+\be\sum_{i\in\hsi}\phi\lp(\hu[i]\rp)+\be\sum_{i\in\hsi^c}\phi\lp(v[i]\rp)\nn\\
\Big[\mb{\rm by~\eq{rer}}\Big]~~~~
    &=&\|A(\hu+v)-d\|^2+\be\sum_{i\in\hsi}\phi\lp(\hu[i]\rp)\nn\\
\Big[\mb{\rm by~\eq{yui}}\Big]~~~~&\geq&\|A\hu-d\|^2+\be\sum_{i\in\hsi}\phi\lp(\hu[i]\rp)\nn\\
\Big[\mb{\rm by~\eq{fds}}\Big]~~~~&=&\Fd(\hu)
~.\lab{4v}\eeqnn
Combining the obtained implication with~\lem{tl}(ii)  shows that
\[\Fd(\hu+v)\geq\Fd(\hu)~~~\all v\in \Bm_\infty(0,\rho)~.\]

If $\hu=0$, this is a (local) minimizer of $\Fd$ by \lem{uz}.
\endproof

Many authors mention that initialization is paramount for the success of approximate algorithms minimizing $\Fd$.
In view of \prop{cp}, if one already has a  well-elaborated
initialization, it could be enough to solve the relevant problem \cpa.

The  statement reciprocal to \prop{cp} is obvious but it helps the presentation.

\BL\lab{ob} For $d\in\RM$ and $\be>0$, let $\Fd$ have a (local) minimum at $\hu$.
Then $\hu$ solves \cphsi for $\hsi\deq\si(\hu)$.
\EL

\proof
Let $\hu$ be a (local) minimizer of $\Fd$. Denote $\hsi\deq\si(\hu)$.
Then $\hu$ solves the problem
\[\min_{u\in\RN}\Big\{\|Au-d\|^2+\be\#\hsi\Big\}~~\mb{subject to } u[i]=0~~\all i\in \hsi^c.\]
Since $\#\hsi$ is a constant, $\hu$ solves \cphsi.
\endproof

\BR ~~   By \prop{cp} and \lem{ob},
solving \cpa for some $\vt\subseteq\IN$ \\
is {\em equivalent} to finding a (local) minimizer of~$\Fd$.
\lab{x1}\ER

This equivalence underlies most of the theory developed in this work.

\BC\lab{nor}
For $d\in\RM$ and $\be>0$, let $\hu$ be a (local) minimizer of $\Fd$. Set $\hsi\deq\si(\hu)$. Then
\beq\lab{ee}\hu=Z_{\hsi}(\uhs)~,~~~\mb{\rm where $~\uhs~$ satisfies}~~~\Ahs^T\Ahs\uhs=\Ahs^Td~.
\eeq
Conversely, if $\hu\in\RN$ satisfies \eq{ee} for  $\hsi=\si(\hu)$, then $\hu$ is a (local) minimizer of $\Fd$.
\EC
\proof By \lem{ob}, $\hu$ solves \cphsi. The equation in~\eq{ee} follows directly from
\eq{eu}. The last claim is a straightforward consequence of \eq{eu} and \prop{cp}.
\endproof

\BR \lab{x2} \rm
Equation \eq{ee} shows that a (local) minimizer $\hu$ of $\Fd$ follows a {\em pseudo}-hard thresholding
scheme\footnote{In a Bayesian setting,
the quadratic data fidelity term in~$\Fd$ models data
corrupted with Gaussian i.i.d. noise.}:
the nonzero part $\uhs$ of $\hu$ is the
least squares solution with respect to the submatrix $\Ahs$ and the whole data vector $d$ is involved in its computation.
Unlike the hard thresholding scheme in \cite{Donoho94},
unsignificant or purely noisy data entries can hardly be discarded from  $\hu$
and they threaten to pollute its nonzero part $\uhs$. See also \rem{ifn}.
\ER

Noisy data $d$ should degrade $\uhs$ and this effect is stronger if
$\Ahs^T\Ahs$ is ill-conditioned \cite{Demoment89}.
The quality of the outcome critically depends on the selected (local) minimizer
and on the pertinence of $A$. 

It may be interesting to evoke another consequence of \prop{cp}:
\BR \lab{x3}
Given $d\in\RM$, for any $\vt\subseteq\IN$, $\Fd$ has a (local) minimizer $\hu$ defined by \eq{ee} and obeying
$\si(\hu)\subseteq\vt$.
\ER

\section{The strict minimizers of $\bm{\Fd}$}\lab{sSM}~~

{\em  We remind, yet again, that no special assumptions on $A\in\RMN$ are adopted.}

Strict minimizers of an objective function enable unambiguous solutions of inverse problems.
The definition below is useful in characterizing the strict minimizers of $\mathcal{F}_d$.

\BD\lab{om} Given a matrix $A\in\RMN$,
for any $r\in\IM$ we define $\Om_r$ as the subset of {\em all} $r$-length supports
that correspond to full column rank $\m\x r$ submatrices of $A$, i.e.,
\beq\mb{\framebox{$\disp{~
\Om_r=\Big\{~\vt\subset\IN~:~~\#\vt=r=\rank\lp(\Avt\rp)~\Big\}~.}$}}
\nn\lab{sm}\eeq
Set $\Om_0=\void$ and define as well
\beq\mb{\framebox{$\disp{~
\Om\deq\bigcup_{r=0}^{\m-1}\Om_r~~\mb{\rm and}~~\Omm\deq\Om\cup\Om_\m~.}$}}\nn\lab{si}\eeq
\ED
\defi{om} shows that for any $r\in\IM$,
\[\rank(A)=r\geq1~~~\Leftrightarrow~~~\Om_r\neq\void~~\mb{\rm and}~~\Om_t=\void~~\all t\geq r+1~.\]

\subsection{How to recognize a strict minimizer of $\bm{\Fd}$?}\lab{rec}

The theorem below gives an exhaustive answer to this question.
\BT\lab{ra}
Given $d\in\RM$ and $\be>0$, let $\hu$ be a (local) minimizer of $\Fd$.
Define
\[\hsi\deq \si(\hu)~.\]
The following statements are {\em equivalent}:
\bit
\item[\rm (i)] The (local) minimum that $\Fd$ has at $\hu$ is {\em strict};
\item[\rm (ii)] \framebox{$\disp{\rank\lp(\Ahs\rp)=\#\hsi}$}~;
\item[\rm(iii)] $\hsi\in\Omm$~.
\eit
If $\hu$ is a strict (local) minimizer of $\Fd$, then it reads
\beq \hu=Z_{\hsi}\lp(\uhs\rp)~~~\mb{for}~~~
\uhs=\lp(\Ahs^T\Ahs\rp)^{-1}\Ahs^Td~\lab{hu}
\eeq
and satisfies $\#\hsi= \|\hu\|_0\leq\m~.$
\ET

\proof We break the proof into four parts.

\paragraph{\rm[(i)$\Rightarrow$(ii)]}~
We recall that by the rank-nullity theorem \cite{Golub96,Meyer00}
\beq\dim\ker\lp(\Ahs\rp)=\#\hsi-\rank\lp(\Ahs\rp)~.\lab{rn}\eeq
Let\footnote{This part can alternatively be proven using \rem{ref}.}
 $\hu\neq0$ be a {\em strict} (local) minimizer of $\Fd$.
Assume that (ii) fails. Then \eq{rn} implies that
\beq \dim\ker\lp(\Ahs\rp)\geq1~.\lab{qa}\eeq
By \lem{ob}, $\hu$ solves \cphsi.
Let $\rho$ read as in~\eq{rho} and let $\Km_{\hsi}$ be defined according to~\eq{km}.
Noticing that
\beq v\in\Km_{\hsi},~~\hsi\neq\void~~~\Rightarrow~~~Av=\Ahs v_{\hsi}~,\lab{ui}
\eeq
\lem{tl}(i) shows that
\beqnn \lp\{\barr{c}v\in \Bm_\infty(0,\rho)\cap\Km_{\hsi}~,\\~~\\
v_{\hsi}\in\ker\lp(\Ahs\rp)\earr\rp.~~~\Rightarrow~~~
\Fd(\hu+v)&=&\|\Ahs\lp(\uhs+v_{\hsi}\rp)-d\|^2+\be\sum_{i\in\hsi}\phi\lp(\hu[i]+v[i]\rp)\\
\Big[\mb{\rm by \lem{tl}(i)}\Big]~~~~
&=&\|\Ahs\uhs-d\|^2+\be\sum_{i\in\hsi}\phi\lp(\hu[i]\rp)+\be\sum_{i\in\hsi^c}\phi\lp(v[i]\rp)\\
\Big[\mb{\rm by \eq{rer}~}\Big]~~~~&=&\|\Ahs\uhs-d\|^2+\be\sum_{i\in\hsi}\phi\lp(\hu[i]\rp)\\
\Big[\mb{\rm by \eq{fds}~}\Big]~~~~&=&\Fd(\hu)~,
\eeqnn
i.e., that $\hu$ is not a strict minimizer, which contradicts (i).
Hence the assumption in~\eq{qa} is false.
Therefore (ii) holds true.

If $\hu=0$, then $\hsi=\void$; hence $\Ahs\in\RR^{\m\x0}$ and $\rank\lp(\Ahs\rp)=0=\#\hsi$ according to \eq{cw}.

\paragraph{\rm$\mathrm{[(ii)\Rightarrow(i)]}$}  Let $\hu$ be a minimizer of $\Fd$ that satisfies (ii).
To have $\#\hsi=0$ is equivalent to $\hu=0$.
By \lem{uz}, $\hu$ is a strict minimizer.
Focus on $\#\hsi\geq1$.
Since $\rank\lp(\Ahs\rp)=\#\hsi\leq\m$ and problem~\cpar in~\eq{cpr} is strictly convex for $\vt=\hsi$,
its unique solution $\uhs$ satisfies
\[v\in\RR^{\#\hsi}\sm\{0\}~~~\Rightarrow~~~\|\Ahs\lp(\uhs+v\rp)-d\|^2>\|\Ahs\uhs-d\|^2~.\]
Using \eq{ui}, this is equivalent to
\beq v\in\Km_{\hsi}\sm\{0\}~~~\Rightarrow~~~\|A(\hu+v)-d\|^2=\|\Ahs\lp(\uhs+v_{\hsi}\rp)-d\|^2>\|\Ahs\uhs-d\|^2
=\|A\hu-d\|^2~.\lab{$}\eeq
 \lem{tl}(i), along with \eq{rer}, yields
\beqnn v\in\Bm_\infty(0,\rho)\cap\Km_{\hsi}\sm\z~~~\Rightarrow~~~\Fd(\hu+v)
&=&\|A(u+v)-d\|^2+\be\sum_{i\in\hsi}\phi\lp(\hu[i]\rp)\\
\Big[\mb{\rm by~\eq{$}}\Big]~~~~~~~~~&>&\|A\hu-d\|^2+\be\sum_{i\in\hsi}\phi\lp(\hu[i]\rp)\\
&=&\Fd(\hu)~.
\eeqnn
Since $\#\hsi\leq\m\leq\n-1$, we have $\hsi^c\neq\void$.
So \lem{tl}(ii) tells us that
\[v\in\Bm_\infty(0,\rho)\sm\Km_{\hsi}~~~\Rightarrow~~~\Fd(\hu+v)>\Fd(\hu)~.\]
Combining the last two implications proves (i).

\paragraph{\rm$\mathrm{[(ii)\Rightarrow(iii)]}$}
Comparing (iii) with Definitions \ref{sp} and \ref{om} proves the claim.

\paragraph{\rm[Equation~\eq{hu}]} The proof
follows from equation~\eq{ee} in \cor{nor}  where\footnote{For
$\hsi=\void$, \eq{zo} and~\eq{cw}$(a)$ show that~\eq{hu} yields  $\hu=0$.}
$\Ahs^T\Ahs$ is invertible.
\endproof

{\em \thm{ra} provides a simple rule enabling one to verify whether or not a numerical
scheme has reached a strict (local) minimizer of $\Fd$.}

The notations $\Om_r$, $\Om$ and $\Omm$ are frequently used in this paper.
Their interpretation is  obvious in the light of \thm{ra}.
For any $d\in\RM$ and for all $\be>0$, the set $\Omm$ is composed of
the supports of all possible strict (local) minimizers of $\Fd$, while $\Om$ is
is the subset of those that are $(\m-1)$-sparse.

An easy and useful corollary is presented next.

\BC\lab{blg}
Let $d\in\RM$. Given an arbitrary $\vt\in\Omm$, let $\hu$ solve \cpa.
Then
\bit
\item[\rm(i)] $\hu$ reads as
\beq\hu=Z_\vt\lp(\hvt\rp)~,~~\mb{where}~~~\hvt=\lp(\Avt^T\Avt\rp)^{-1}\Avt^Td~,\lab{bla}\eeq
and obeys $\hsi\deq\si(\hu)\subseteq\vt$ and $\hsi\in\Omm$~;
\item[\rm(ii)] for any $\be>0$, $\hu$ is a {\em strict} (local) minimizer of $\Fd$;
\item[\rm(iii)] $\hu$ solves \cphsi.
\eit
\EC

\proof
Using~\eq{eu}, $\hu$ fulfills (i) since $\Avt^T\Avt$
is invertible and $\si(\hu)\subseteq\vt$ by the constraint in \cpa.
If $\hsi=\void$, (ii) follows from~\lem{uz}.
For $\#\hsi\geq1$, $\Ahs$ is an $\m\x\#\hsi$ submatrix of $\Avt$.
Since $\rank\lp(\Avt\rp)=\#\vt$, we have $\rank\lp(\Ahs\rp)=\#\hsi$ and so $\hsi\in\Omm$.
By~\prop{cp} $\hu$ is a (local) minimizer of $\Fd$, and~\thm{ra} leads to (ii).
 \lem{ob} and \cor{blg}(ii) yield (iii).
\endproof

\BR \lab{x4} One can easily compute a strict (local) minimizer $\hu$ of $\Fd$ without knowing
the value of the regularization parameter $\be$.
Just data $d$ and an $\vt\in\Omm$ are needed.
\ER

This consequence of~\cor{blg} might be striking.

\ms
Clearly,  the support $\sigma(\bar u)$ of a nonstrict local minimizer $\bu$ of $\Fd$ contains some subsupports
yielding strict (local) minimizers of $\Fd$. It is easy to see that
among them, there is $\hsi\subsetneqq\si(\bu)$ such that the corresponding $\hu$ given by \eq{hu}
{\em strictly} decreases the value of  $\Fd$; i.e., $\Fd(\hu)<\Fd(\bu).$

\subsection{Every strict (local) minimizer of $\bm{\Fd}$ is linear in $\bm{d}$}\lab{els}
Here we explore the behavior of the strict (local) minimizers of $\Fd$ with respect to variations of~$d$.
An interesting sequel of \thm{ra} is presented in the following corollary.

\BC\lab{sk}
For $d\in\RM$ and $\be>0$, let $\hu$ be a (local) minimizer of $\Fd$ satisfying
$\hsi\deq\si(\hu)\in\Om~.$
Define
\[\Nm_{\hsi}\deq\ker\lp(\Ahs^T\rp)\subset\RM~.\]
We have $\,\dim\Nm_{\hsi}=\m-\#\hsi\geq1$ and
\[{d^{\,\prime}}\in\Nm_{\hsi}~~~\Rightarrow~~~\F_{d+{d^{\,\prime}}}~~\mb{has a strict (local) minimum at}~~\hu~.\]
\EC

\proof
Since $\hsi\in\Om$, the minimizer $\hu$ is strict by \thm{ra}.
By ${d^{\,\prime}}\in\ker\lp(\Ahs^T\rp)$ we find $\Ahs^T(d+{d^{\,\prime}})=\Ahs^Td$ for any $\,{d^{\,\prime}}\in\Nm_{\hsi}$.
Inserting this into~\eq{hu} in \thm{ra} yields the result.
\endproof

All data located in the vector subspace $\Nm_{\hsi}\supsetneqq\{0\}$ yield the same strict
(local) minimizer $\hu$.

\begin{remark}  \lab{ifn} {\rm If data contain noise $n$, it can be decomposed in a unique way as
$n=n_{\Nm_{\hsi}}+n_{\Nm_{\hsi}^\bot}$ where $n_{\Nm_{\hsi}}\in\Nm_{\hsi}$
and $n_{\Nm_{\hsi}^\bot}\in\Nm_{\hsi}^\bot$.
The component $n_{\Nm_{\hsi}}$ is removed (\cor{sk}), while $n_{\Nm_{\hsi}^\bot}$
is transformed according to~\eq{hu} and added to $\uhs$. }
\end{remark}

We shall use the following definition.

\BD Let ${\O}\subseteq\RM$ be an open domain.
We say that $\U:{\O}\to\RN$ is a {\em local minimizer function} for the
family of objectives $\mathfrak{F}_{\O}\deq\{\Fd~:~d\in {\O}\}$ if for any $d\in {\O}$,
the function $\Fd$ reaches a {\em strict} (local)  minimum at $\U(d)$.
\lab{mf}
\ED

\cor{blg} shows that for any $d\in\RM$, each strict (local) minimizer of $\Fd$ is entirely
described by an $\vt\in\Omm$ via equation~\eq{bla} in the same corollary.
Consequently, a local minimizer function $\U$ is associated with every $\vt\in\Omm$.

\BL\lab{ct}
For some arbitrarily fixed $\vt\in\Omm$ and $\be>0$, the family of functions $\mathfrak{F}_{\RM}$
has a {\em linear} (local) minimizer function  $\U:\RM\to\RN$ that reads as
\beq \all d\in\RM,~~~\U(d)=Z_\vt\lp(U_\vt \,d\rp)~,~~
\mb{where}~~
    U_\vt=\lp(\Avt^T\Avt\rp)^{-1}\Avt^T ~\in\RR^{\#\vt\x\m}~.
\lab{bu}\eeq
\EL

\proof
The function $\U$ in~\eq{bu} is linear in $d$.
From~\cor{blg}, for any $\be>0$ and for any $d\in\RM$, $\Fd$ has a strict (local) minimum
at $\hu=\U(d)$.
Hence $\U$ fits \defi{mf}.
\endproof

{\em Thus, even if $\Fd$ has many strict (local) minimizers, each is linear in $d$.}

Next we exhibit a {\em closed negligible} subset of $\RM$, associated with a nonempty $\vt\in\Omm$,
whose elements are data $d$ leading to $\|\U(d)\|_0<\#\vt$.

\BL\lab{dom}
For any $\vt\in\Omm$, define the subset $\Dm_\vt\subset\RM$ by
\beq\Dm_\vt\deq\bigcup_{i=1}^{\#\vt}
\Big\{~g\in\RM~:~\lp\<e_i\,,\,\lp(\Avt^T\Avt\rp)^{-1}\Avt^T \;g\rp\>=0 ~\Big\}~.\lab{can}\eeq
Then $\Dm_\vt$ is closed in $\RM$ and $\LL^\m\lp(\Dm_\vt\rp)=0$.
\EL

\proof  If $\vt=\void$ then $\Dm_\vt=\void$, hence the claim.
Let $\#\vt\geq1$.
For some $i\in\II_{\#\vt}$, set
\beqnn\Dm&\deq&\Big\{~g\in\RM~:~\lp\<e_i\,,\,\lp(\Avt^T\Avt\rp)^{-1}\Avt^T \;g\rp\>=0 ~\Big\}\\&=&
\Big\{~g\in\RM~:~\Big\<\Avt\lp(\Avt^T \Avt\rp)^{-1}e_i,\;g\,\Big\>=0\Big\}~.
\eeqnn
Since $\rank\!\lp(\Avt\lp(\Avt^T \Avt\rp)^{-1}\rp)\!=\!\#\vt$,
$\ker\!\lp(\Avt\lp(\Avt^T \Avt\rp)^{-1}\rp)\!=\!\z$.
Hence  $A_\omega\left(A_\omega^T A_\omega\right)^{-1}\!e_i\neq0$.
Therefore $\Dm$ is a vector subspace of $\RM$ of dimension $\m-1$
and so $\LL^\m\lp(\Dm\rp)=0$.
The conclusion follows from the fact that $\Dm_\vt$ in~\eq{can} is the  union of
$\#\vt$ subsets like $\Dm$ (see, e.g., \cite{Rudin76,Evans92}).
\endproof

\BP\lab{ctp} For some arbitrarily fixed $\vt\in\Omm$ and $\be>0$,
let $\,\U:\RM\to\RN$ be the relevant (local)
minimizer function for $\mathfrak{F}_{\RM}$  as given in~\eq{bu} (\lem{ct}).
Let $\Dm_\vt$ read as in \eq{can}.
Then the function $d\mapsto\Fd\big(\U(d)\big)$ is $\C^\infty$ on $\RM\sm\Dm_\vt$ and
\beq d\in\RM\sm\Dm_\vt~~~\Rightarrow~~~
\si\lp(\U(d)\rp)=\vt~,\nn\lab{mp}\eeq
where the set $\RM\sm\Dm_\vt$ contains an open and dense subset of $\RM$.
\EP

\proof The statement about $\RM\sm\Dm_\vt$ is a direct consequence of \lem{dom}.

If $\vt=\void$, then $\U(d)=0$  for all $d\in\RM$, so all claims in the proposition are trivial.
Consider that $\#\vt\geq1$.
For any $i\in\II_{\#\vt}$, the $\vt[i]$th component of $\,\U(d)$ reads as (see \lem{ct})
\[\U_{\vt[i]}(d)=\lp\<e_i\,,\,\lp(\Avt^T\Avt\rp)^{-1}\Avt^T \;d\rp\>~.\]
The definition of $\Dm_\vt$ shows that
\[d\in\RM\sm\Dm_\vt~~\mb{\rm and}~~i\in\II_{\#\vt}~~~\Rightarrow~~~
\U_{\vt[i]}(d)\neq0~,\]
whereas $\U_i(d)=0$ for all $i\in\vt^c$.
Consequently,
\[\vt\in\Omm~~\mb{\rm and}~~d\in\RM\sm\Dm_\vt~~~\Rightarrow~~~\si\lp(\U(d)\rp)=\vt~.\]
Then $\|\U(d)\|_0$ is constant on $\RM\sm\Dm_\vt$ and
\[\vt\in\Omm~~\mb{\rm and}~~d\in\RM\sm\Dm_\vt~~~\Rightarrow~~~\Fd(\U(d))=\big\|A\U(d)-d\big\|^2+\be\#\vt~.\]
We infer from \eq{bu} that $d\mapsto\big\|A\U(d)-d\big\|^2$ is
$\C^\infty$ on $\RM$, so $d\mapsto\Fd\big(\U(d)\big)$ is $\C^\infty$ on $\RM\sm\Dm_\vt$.
\endproof

{\em A generic property is that a local minimizer function corresponding to $\Fd$ produces solutions sharing the same support.}
The application $d\mapsto\Fd\big(\U(d)\big)$ is discontinuous on the {\em closed negligible} subset
$\Dm_\vt$, where the support of $\U(d)$ is shrunk.

\subsection{Strict minimizers with an $\m$-length support} \lab{smm}

Here we explain why minimizers with an $\m$-length support are useless in general.

\BP\lab{bad}
Let $\rank(A)=\m$, $\be>0$ and for $d\in\RM$ set
\[\Um_\m\deq\big\{\hu\in\RN~:~ \hu~~\mb{is a strict (local) minimizer of ~$\Fd$~ meeting}~~\|\hu\|_0=\m~\big\}~.\]
Put
\beq\Qm_\m\deq\bigcup_{\vt\in\Om_\m}~\bigcup_{i\in\IM}
\big\{g\in\RM~:~\lp\<e_i\,,\,\Avt^{-1}g\rp\>=0\,\big\}~.
\lab{qm}\eeq
Then $\RM\sm\Qm_\m$ contains a dense open subset of $\RM$ and
\[d\in\RM\sm\Qm_\m\quad\Rightarrow\quad\#\Um_\m=\#\Om_\m\quad\mb{\rm and}\quad\Fd(\hu)=\be\m~~\all\hu\in\Um_\m~.\]
\EP

\proof
Using the notation in~\eq{can}, $\Qm_\m$ reads as
\[\Qm_\m=\bigcup_{\vt\in\Om_\m}\Dm_\vt~.\]
The claim on $\RM\sm\Qm_\m$ follows from~\lem{dom}.
Since $\rank(A)=\m$, we have $\#\Om_\m\geq1$.

Consider that $d\in\RM\sm\Qm_\m$.
By \prop{ctp}
\[d\in\RM\sm\Qm_\m\quad\mb{\rm and}\quad\vt\in\Om_\m\quad\Rightarrow\quad\Fd~
\mb{has a strict (local) minimizer ~$\hu$~ obeying}~\si(\hu)=\vt~.\]
Hence $\hu\in\Um_\m$.
Therefore, we have a mapping $b:\Om_\m\to\Um_\m$ such that $\hu=b(\vt)\in\Um_\m$.
Using \lem{ob} and \cor{blg}, it reads as
\[b(\vt)=Z_\vt(\Avt^{-1}d)~.\]
For $(\vt,\vs)\in\Om_\m\x\Om_\m$ with $\vs\neq\vt$ one obtains $\hu=b(\vt)\in\Um_\m$,
$\bu=b(\vs)\in\Um_\m$ and $\hu\neq\bu$, hence $b$ is one-to-one.
Conversely,
for any $\hu\in\Um_\m$ there is $\vt\in\Om_\m$ such that $\hu=b(\vt)$ and $\si(\hu)=\vt$ (because $d\not\in\Qm_\m$).
It follows that $b$ maps $\Om_\m$ onto $\Um_\m$.
Therefore, $\Om_\m$ are $\Um_\m$ in one-to-one correspondence, i.e. $\#\Om_\m=\#\Um_\m$.

Last, it is clear that $\vt\in\Om_\m$ and $d\not\in\Qm_\m$ lead to
$\|A\hu-d\|^2=0$ and  $\Fd(\hu)=\be\m$.
\endproof

{\em For any $\be>0$, a generic property of
$\Fd$ is that it has $\#\Om_\m$
strict minimizers $\hu$ obeying $\|\hu\|_0=\m$ and  $\Fd(\hu)=\be\m$.
It is hard to discriminate between all these  minimizers.
Hence the interest in minimizers with supports located in $\Om$, i.e.,
strict $(\m-1)$-sparse minimizers of $\Fd$.}

\section{On the global minimizers of $\bm{\Fd}$}\lab{sgm}

The next proposition gives a {\em necessary condition} for a global minimizer of $\Fd$.
It follows directly from~\cite[Proposition 3.4]{Nikolova04a}
where\footnote{
Just set $g_i=e_i$, $P=\Id_\m$
and $H=\Id_\n$ in \cite[Proposition 3.4]{Nikolova04a}.}
the regularization term is $\|Du\|_0$ for a  full row rank matrix $D$.
For  $\Fd$ in~\eq{fd} with $\|a_i\|_2=1$, $\all i\in\IN$,
a simpler condition was derived later in \cite[Theorem 12]{Tropp06}, using different tools.
For completeness, the proof for a general $A$ is outlined in Appendix~\ref{pou}.

\BP\lab{ou}
For $d\in\RM$ and $\be>0$, let $\Fd$ have a global minimum at  $\hu$. Then\footnote{Recall that
$a_i\neq0$ for all $i\in\IN$  by \eq{ao}  and that $\|\cdot\|=\|\cdot\|_2$.}
\beq i\in\si(\hu)~~~\Rightarrow~~~\big|\,\hu[i]\,\big|\geq \frac{\sqrt{\be}}{\|a_i\|}~.
\lab{rt}\eeq
\EP

Observe that the lower bound on
$\disp{\lp\{\big|\,\hu[i]\,\big|~:~i\in\si(\hu)\rp\}}$ given in
\eq{rt} is independent of~$d$.
This means that in general, \eq{rt} provides a pessimistic bound.

 The proof of the statement shows that~\eq{rt} is met also by all (local) minimizers of $\Fd$~satisfying
$$\Fd(\hu)\leq\Fd(\hu+\rho e_i)~~~\all\rho\in\RR,~~\all i\in\IN~.$$

\subsection{The global minimizers of $\bm{\Fd}$ are strict}\lab{sGM}~~

\begin{remark}  \lab{nu} \rm
Let $d\in\RM$ and $\be>\|d\|^2$.
Then $\Fd$ has a strict global minimum at $\hu=0$. {\rm Indeed,}
\[v\neq 0~~\Rightarrow~~\|v\|_0\geq1~~\Rightarrow~~
\Fd(0)=\|d\|^2<\be\leq\|Av-d\|^2+\be\|v\|_0~.\]
\end{remark}

For least-squares regularized with a more regular $\phi$, one usually gets $\hu=0$
asymptotically as $\be\to+\infty$ but $\hu\neq0$ for finite values of $\be$.
This does not hold for $\Fd$ by \rem{nu}.

Some theoretical results on the global minimizers of $\Fd$ have been obtained
\cite{Nikolova04a,Haupt06,Tropp06,Davies08}.
Surprisingly, the question about the {\em existence of global minimizers of $\Fd$} has never been raised.
We answer this question using the notion of {\em asymptotically level stable functions}
introduced by Auslender~\cite{Auslender00} in 2000.
As usual,
\[\lev(\Fd,\la)\deq\{v\in\RN~:~\Fd(v)\leq\la\}\qu\mb{\rm for}\qu \la>\inf \Fd~.\]
The following definition is taken from \cite[p.~94]{Auslender03}.

\BD\lab{als}
Let $\Fd:\RN\to\RR\cup\{+\infty\}$ be lower semicontinuous and proper.
Then $\Fd$ is said to be {\em asymptotically level stable} if for each $\rho>0$, each bounded
sequence $\{\la_k\}\in\RR$ and each sequence $\{v_k\}\in\RN$ satisfying
\beq v_k\in\lev(\Fd,\la_k),\qu\|v_k\|\to+\infty, \qu v_k\;\|v_k\|^{-1}\to \bv\in\ker\lp((\Fd)_\infty\rp)~,
\lab{aal}\eeq
where $(\Fd)_\infty$ denotes the asymptotic (or recession) function of $\Fd$, there exists $k_0$ such that
\beq v_k-\rho\bv\in\lev(\Fd,\la_k)\qu\all k\geq k_0~.\lab{lsa}\eeq
\ED

One can note that a coercive function is  asymptotically level stable, since \eq{aal} is empty.
We prove that our discontinuous noncoercive objective $\Fd$ is asymptotically level stable as well.

\BP\lab{pal} Let $\Fd:\RN\to\RR$ be of the form~\eq{fd}. Then $\ker\lp((\Fd)_\infty\rp)=\ker(A)$ and
$\Fd$ is asymptotically level stable.
\EP

The proof is outlined in Appendix~\ref{palp}.

\BT\lab{gs}Let $d\in\RM$ and $\be>0$. Then
\bit
\item[\rm(i)]
the set of the global minimizers of $\Fd$
\beq\disp{\h U\deq \lp\{\hu\in\RN~:~\hu=\min_{u\in\RN}\Fd(u)\rp\}}\lab{ops}\eeq is nonempty;
\item[\rm(ii)] every $\hu\in \h U$ is a {\em strict} minimizer of $\Fd$, i.e.,
\[\si(\hu)\in\Omm~,\]
hence $\|\hu\|_0\leq\m$.
\eit
\ET

\proof
For claim (i), we use the following statement\footnote{This
result was originally exhibited in \cite{Baiocchi98} (without the
notion of asymptotically level stable functions).},
whose proof can be found in the monograph by Ausleneder and Teboulle~\cite{Auslender03}:

\noi{\cite[Corollary 3.4.2]{Auslender03}} {\it
Let $\Fd:\RN\to\RR\cup\{+\infty\}$ be asymptotically level stable
with $\inf\Fd>-\infty$.
Then the optimal set $\h U$---as given in \eq{ops}---is nonempty .}

\noi From \prop{pal}, $\Fd$ is asymptotically level stable and $\inf\Fd\geq0$ from~\eq{fd}.
Hence $\h U\neq\void$.

(ii).~ Let $\hu$ be a global minimizer of $\Fd$. Set $\hsi=\si(\hu)$.

If $\hu=0$, (ii) follows from \lem{uz}.
Suppose that the global minimizer $\hat u\neq0$ of $\mathcal{F}_d$ is {\em nonstrict}.
Then \thm{ra}(ii) fails to hold; i.e.,
\beq\dim\ker\lp(\Ahs\rp)\geq1~. \lab{sit}\eeq
Choose $v_{\hsi}\in \ker\lp(\Ahs\rp)\sm\{0\}$ and set $v=Z_{\hsi}\lp(v_{\hsi}\rp)$.
Select an $i\in\hsi$ obeying $v[i]\neq0$.
Define $\wt u$ by
\beq \wt u\deq \hu -\hu [i]\frac{v}{v[i]}~.\eeq
We have $\wt u[i]=0$ and $\hu [i]\neq0~.$
 Set $\wt\si\deq\si\lp(\wt u\rp)$.
Then
\beq\wt\si\subsetneqq\hsi~~~\mb{hence}~~~\#\wt\si\leq\#\hsi-1~.\eeq
From $\disp{v_{\hsi}\frac{\hu [i]}{v[i]}\in\ker\lp(\Ahs \rp)}$, using~\eq{sd}
and \rem{ref} shows that\footnote{In detail we have
$A\hu =\Ahs \uhs =\Ahs \lp(\uhs -v_{\hsi}\frac{\hu [i]}{v[i]}\rp)=\Ahs \wt u_{\hsi}=
A_{\wt\si}\wt u_{\wt\si}=A\wt u~.$}
$A\hu =\Ahs \uhs =A_{\wt\si}\wt u_{\wt\si}=A\wt u$.
Then
\beqn\Fd\lp(\wt u\rp)&=&\|A\wt u-d\|^2+\be\#\wt\si\nn\\
&\leq&\Fd\lp(\hu \rp)-\be=\|A\hu -d\|^2+\be\lp(\#\hsi-1\rp)~.\nn\lab{tr}\eeqn
It follows that $\hu$ is not a global minimizer, hence \eq{sit} is false.
Therefore $\rank\lp(\Ahs\rp)=\#\hsi$ and $\hu$ is a strict minimizer of $\Fd$ (\thm{ra}).
\endproof

One can note that if $\rank(A)=\m$, any global minimizer $\hu$ of $\Fd$ obeys
$\,\Fd(\hu)\leq\be\m~.$

  According to \thm{gs}, {\em the global minimizers of $\Fd$ are strict and their number is finite: this is
a nice property that fails for many convex nonsmooth optimization problems.}

\subsection{$\bm\k$-sparse global minimizers for $\bm{\k\leq\m-1}$}

In order to simplify the presentation, in what follows we consider that
\[\mb{\framebox{~$\rank(A)=\m<\n~.$}}\]

Since $\Fd$ has a large number (typically equal to $\#\Om_\m$)
 of strict minimizers with $\|\hu\|_0=\m$
yielding the same value $\Fd(\hu)=\be\m$ (see \prop{bad} and the comments given after its proof),
it is important to be sure that the global minimizers of $\Fd$ are $(\m-1)$-sparse.

We introduce a notation which is used in the rest of this paper.
For any $\k\in\Im$, put
\beq \mb{\framebox{$\disp{~~\Ombk\deq\bigcup_{r=0}^\k\Om_r~}$~}}\lab{ga}\eeq
where $\Omega_r$ was set up in~\defi{om}.
\thm{ra} gives a clear meaning of the sets\footnote{Clearly, $\Omb_{\m-1}=\Om~.$}~$\Ombk$.
{\em For any $d\in\RM$ and any $\be>0$, for any fixed $\k\in\Im$, the set $\Ombk$ is composed
of the supports of all possible $\k$-sparse strict (local) minimizers of $\Fd$.}

The next propositions checks the existence of $\be>0$ ensuring that all
the global minimizers of $\Fd$ are $\k$-sparse, for some $\k\in\Im$.

\BP\lab{ie} Let $d\in\RM$.
 For any $\k\in\Im$, there exists $\be_\k\geq0$ such that if $\be>\be_\k$, then
each global minimizer $\hu$ of $\Fd$ satisfies
\beq\|\hu\|_0\leq\k~~~\mb{\rm and}~~~\si(\hu)\in\Ombk~ .\lab{spo}\eeq
One can choose $\be_\k=\|A\tu-d\|^2$ where $\tu$ solves \cpa for some $\vt\in\Om_{\k}$.
\EP

The proof is given in Appendix~\ref{pie}.
The value of $\be_\k$ in the statement is easy to compute, but in general it
is not sharp\footnote{
For $\be\gtrapprox\be_\k$
the global minimizers of $\Fd$ might be $k$-sparse for $k\ll\k$.
A sharper $\be_\k$ can be obtained~as
$\be_\k=\min_{\vt\in\Om_\k}\lp\{\|A\tu-d\|^2~:~\tu~~\mb{solve \cpa for}~~ \vt\in\Om_\k\rp\}~.$
}.

\section{Uniqueness of the global minimizer of $\bm{\Fd}$}\lab{kb}

The presentation is simplified using the notation introduced next.
Given a matrix $A\in\RMN$,
with any $\vt\in\Om$ (see \defi{om}), we associate the $\m\x\m$ matrix $\Pvt$ that yields the orthogonal
projection\footnote{
If $\vt=\void$, we have $\Avt\in\RR^{\m\x0}$ and so $\Pi_\omega$ is an  $\m\x\m$ null matrix.}
onto the subspace spanned by the columns of $\Avt$:
\beq\mb{\framebox{$\disp{~\Pvt\deq\Avt\lp(\Avt^T\Avt\rp)^{-1}\Avt^T~.}$}}\lab{pro}\eeq
\noi For  $\vt\in\Om$, the projector in \rem{ref} reads
$\Pi_{\ran(\Avt)}= \Pvt$.

Checking whether a global minimizer $\hu$ of $\Fd$ is unique requires us to compare its value
$\Fd(\hu)$ with the values $\Fd(\bu)$ of the concurrent strict minimizers~$\bu$.
Let $\hu$ be an $(\m-1)$-sparse strict (local) minimizer of $\Fd$. Then $\hsi\deq\si(\hu)\in\Om$.
Using~\rem{ref} shows that
\beqn \Fd(\hu)&=&\|\Ahs\uhs-d\|^2+\beta\#\hsi
=\|\Phs d-d\|^2+\beta\#\hsi\nn\\&=&d^T\lp(\Id_\m-\Phs\rp)d+\beta\#\hsi~.\lab{sx}\eeqn
Let $\bu$ be another $(\m-1)$-sparse strict minimizer of $\Fd$; set $\bsi\deq\si(\bu)$.
Then
\[\Fd(\hu)-\Fd(\bu)=d^T\lp(\Pbs-\Phs\rp)d+\beta(\#\hsi-\#\bsi)~.\]
If both $\hu$ and $\bu\neq\hu$ are  global minimizers of $\Fd$, the previous equality yields
\beq \Fd(\hu)=\Fd(\bu)~~~\Leftrightarrow~~~ d^T\lp(\Pbs-\Phs\rp)d=-\beta(\#\hsi-\#\bsi)~.\lab{hq}\eeq
{\em Equation~\eq{hq} reveals that the uniqueness of the global minimizer of $\Fd$
cannot be guaranteed without suitable assumptions on $A$ and on~$d$.}

\subsection{A generic assumption on $\bm{A}$}\lab{ru}

We adopt an assumption on the matrix $A$ in $\Fd$
in order to restrict the cases when~\eq{hq} takes place for some supports $\hsi\neq\bsi$ obeying $\#\hsi=\#\bsi$.
\BH\lab{aa} The matrix $A\in\RMN$, where $\mathsf{M}<\mathsf{N}$, is such that for some given $\k\in\Im$,
\beq \mbox{\framebox{$ ~ r\in\mathbb{I}_{\mathsf{K}}~~\mbox{\rm and}~~
(\omega,\varpi)\in(\Omega_r\times \Omega_r),~~\omega\neq\varpi~~~\Rightarrow~~~\Pi_{\omega}\neq \Pi_{\varpi}~.$}}\eeq
\EH

Assumption H\ref{aa} means that
the angle (or the gap) between the equidimensional subspaces
$\ran\lp(\Avt\rp)$ and $\ran\lp(\Avs\rp)$ must be nonzero~\cite{Meyer00}.
For instance, if $(i,j)\in\IN\x\IN$ satisfy $i\neq j$, H\ref{aa} implies that $a_i\neq\ka\,a_j$ for any $\ka\in\RR\sm\z$
since $\Pi_{\{i\}}=a_ia_i^T/\|a_i\|^2~.$

Checking whether H\ref{aa} holds for a given matrix $A$ requires a combinatorial search over all possible
couples $(\vs,\vt)\in(\Om_r\x \Om_r)$ satisfying $\vs\neq\vt$, $\all r\in\Ik$.
This is hard to do.
Instead, we wish to know whether or not H\ref{aa} is a {\em practical} limitation.
Using some auxiliary claims, we shall show that H\ref{aa} {\em fails only for a
closed negligible subset of matrices} in the space of all $\m\x\n$ real matrices.

\BL\lab{bb}
Given $r\in\Im$ and $\vs\in\Om_r$, define the following set of submatrices of $A$:
\[\H_\vs=\Big\{\Avt~:~\vt\in\Om_r~~\mb{\rm and}~~\Pvt=\Pvs\Big\}~.\]
Then $\H_\vs$ belongs to an $(r\x r)$-dimensional subspace of the space  of all $\m\x r$ matrices.
\EL

\proof Using  the fact that $\vs\in\Om_r$ and $\vt\in\Om_r$, we have\footnote{Using \eq{pro},
as well as the fact that $\Avt=\Pvt\Avt$, $\all\vt\in\Om_r$, one easily derives \eq{fo} since
\[\Pvt=\Pvs~~\Leftrightarrow\left\{\barr{l}
\Avt\lp(\Avt^T\Avt\rp)^{-1}\Avt^T=\Pvs\\\Avs\lp(\Avs^T\Avs\rp)^{-1}\Avs^T=\Pvt\earr\right.\Rightarrow~
\left\{\barr{l}\Avt=\Pvs\Avt\\\Avs=\Pvt\Avs\earr\right.\Rightarrow~
\left\{\barr{l}\Pvt=\Pvs\Pvt\\\Pvs=\Pvt\Pvs\earr\right.\Rightarrow~\Pvt=\Pvs.\]}
\beq\Pvt=\Pvs~~~\Leftrightarrow~~~\Avt=\Pvs\Avt~~.\lab{fo}\eeq
Therefore $\H_\vs$ equivalently reads
\beq\H_\vs=\Big\{\Avt~:~\vt\in\Om_r~~\mb{\rm and}~~\Avt=\Pvs\Avt\Big\}~.\lab{gr}\eeq
Let $\Avt\in\H_\vs$.
Denote the columns of $\Avt$ by $\t a_i$ for $i\in\Ir$.
Then \eq{gr} yields
\[ \Pvs\t a_i=\t a_i,~~\all i\in\Ir\quad\Rightarrow\quad\t a_i\in\ran(\Avs), ~~\all i\in\Ir~.\]
Hence all $\t a_i$, $i\in\Ir$, live in the $r$-dimensional vector subspace $\ran(\Avs)$.
All the columns of each matrix $\Avt\in\H_\vs$ belong to $\ran(\Avs)$ as well.
It follows that $\H_\vs$
belongs to a (closed) subspace of dimension $r\x r$ in the space of
all $\m\x r$ matrices, where $r\leq\m-1$.
\endproof

More details on the submatrices of $A$ living in  $\H_\vs$ are given next.
\begin{remark}  \lab{P} {\rm The {\em closed negligible} subset $\H_\vs$ in~\lem{bb} is formed from all the
submatrices of $A$ that are column equivalent to $\Avs$
(see \cite[p.~171]{Meyer00}), that is,
\beq\Avt\in\H_\vs~~~\Leftrightarrow~~~\exists\,P\in\RR^{r\x r}~~\mb{such that}~~
\rank(P)=r~~\mb{\rm and}~~\Avt=\Avs P~.\lab{yr}
\eeq
Observe that $P$ has $r^2$ unknowns that must satisfy $\m r$ equations
and that $P$ must be invertible.
It should be quite unlikely that such a matrix $P$ does exist.}
\end{remark}

This remark can help to discern whether or not structured dictionaries satisfy H\ref{aa}.

Next we inspect the set of all matrices $A$ failing assumption H\ref{aa}.
\BL\lab{ug} Consider the set $\H$ formed from $\m\x\n$ real matrices described next:
\[\H\deq\lp\{A\in\RMN~:\exist r\in\Im,~
\exist(\vs,\vt)\in\Om_r\x\Om_r,~\vs\neq\vt~~\mb{\rm and}~~\Pvs=\Pvt\rp\}~.\]
Then $\H$ belongs to a finite union of vector subspaces in $\RMN$ whose Lebesgue measure
in $\mathbb{R}^{\mathsf{M}\times\mathsf{N}}$  is null.
\EL

\proof
Let $A\in\H$. Then there exist at least one integer $r\in\Im$ and at least one pair
$(\vs,\vt)\in\Om_r\x\Om_r$ such that $\vs\neq\vt$  and $\Pvs=\Pvt$.
Using \lem{bb}, $A$ contains (at least) one $\m\x r$ submatrix $\Avs$ belonging to an
$r\x r$ vector subspace in the space of all $\m\x r$ real matrices.
Identifying $A$ with an $\m\n$-length vector, its
entries are included in a vector subspace of $\RR^{\m\n}$ of dimension no larger than $\m\n-1$.
The claim of the lemma is straightforward.
\endproof

We can now clarify assumption H\ref{aa} and show that it is really good.

\BT\lab{usg} Given an arbitrary $\k\in\Im$, consider the set of $\m\x\n$ real matrices below
\[\A_\k\deq\Big\{A\in\RMN~:~~A~~\mb{\rm satisfies~~ H\ref{aa}~~for~~}\k~\Big\}~.\]
Then $\A_\k$ contains an open and dense subset in the space of all $\m\x\n$ real-valued matrices.
\ET

\proof
The complement of $\A_\k$ in the space of all $\m\x\n$ real matrices reads as
\[\A_\k^c=\Big\{A\in\RMN~:~\mb{\rm H\ref{aa} fails for }~A~~\mb{\rm and}~~\k~\Big\}~.\]
It is clear that $\A_\k^c\subset\H~,$ where $\H$ is described in \lem{ug}.
By the same lemma,
$\A_\k^c$ is included in a closed subset of vector subspaces in $\RMN$ whose Lebesgue measure
in $\mathbb{R}^{\mathsf{M}\times\mathsf{N}}$  is null.
Consequently, $\A_\k$ satisfies  the statement of the theorem.
\endproof

{\em For any $\k\in\Im$,
 H\ref{aa} is a {\em generic} property of all $\m\x\n$ real matrices meeting $\m<\n$.}
 This is the meaning of \thm{usg} in terms of \defi{ps}.

We can note that
\[\A_{\k+1}\subseteq\A_\k,~~~\all\k\in\II_{\m-2}~.\]
One can hence presume that H\ref{aa} is weakened as $\k$ decreases.
This issue is illustrated in section~\ref{nm}.

\subsection{A generic assumption on $\bm{d}$}\lab{rs}

A preliminary result is stated next.
\BL\lab{H}
Let $(\vt,\vs)\in\Ombk\x \Ombk$ for $\vt\neq\vs$ and let H\ref{aa} hold for $\k\in\Im$.
Given $\ka\in\RR$, define
\beq\Tm_\ka\deq\{g\in\RM~:~g^T\lp(\Pvt-\Pvs\rp)g=\ka\}~.\nn\lab{setS}\eeq
Then $\Tm_\ka$  is a closed subset of $\RM$ and $\LL^\m\lp(\Tm_\ka\rp)=0$.
\EL

\proof Define $f:\RM\to\RR$ by
$f(g)=g^T\lp(\Pvt-\Pvs\rp)g~.$
Then \beq\Tm_\ka=\{g\in\RM~:~f(g)=\ka\}~.\lab{tyk}\eeq
Using H\ref{aa}, $\Tm_\ka$ is closed in $\RM$.
Set
\[\Qm=\{g\in\RM~:~\nabla f(g)\neq0\}\qu\mb{and}\qu\Qm^c=\RM\sm\Qm.\]

Consider an arbitrary $g\in\Tm_\ka\cap\Qm$. From  H\ref{aa}, $\rank (\nabla f(g))=1$.
For simplicity, assume that
\[\nabla f(g)[\m]=\frac{df(g)}{dg[\m]}\neq0~.\lab{asm}\]
By the implicit functions theorem, there are open
neighborhoods $\O_g\subset\Qm\subset\RM$ and $\V\subset\RR^{\m-1}$ of $g$
and $g_{\Im}$, respectively,
 and a unique $\C^1$-function $h_g:\V\to \RR$ with $\nabla h_g$ bounded, such that
\beq \ga=(\ga_{\Im},\ga[\m])\in\O_g~~\mb{and}~~f(\ga)=\ka\qu\Leftrightarrow\qu
\ga_{\Im}\in\V ~~\mb{and}~~\ga[\m]=h_g(\ga_{\Im})~.\lab{imp}
\eeq
From~\eq{tyk} and \eq{imp} it follows that\footnote{From~\eq{imp}, $\V$ is the restriction of
$\O_g$ to $\RR^{\m-1}$.}
\[\O_g\cap\Tm_k=\psi^g\lp(\O_g\cap(\RR^{\m-1}\x\z)\rp),~\]
where $\psi^g$ is a diffeomorphism on $\O_g $ given by
\[\psi^g_i(\ga)=\ga[i],\qu 1\leq i\leq\m-1\qu\mb{and}\qu\psi^g_\m(\ga)=h_g(\ga_{\Im})+\ga[\m]~.\lab{nap}\]
Since $\LL^\m\lp(\O_g\cap(\RR^{\m-1}\x\z)\rp)=0$ and $\nabla\psi^g$ is bounded on $\O_g$, it
follows from \cite[Lemma 7.25]{Rudin87} that\footnote{The same result follows from
the change-of-variables theorem for the Lebesgue integral, see e.g. \cite{Rudin87}.}
 $\LL^\m\lp(\V_g\cap\Tm_k\rp)=0$.
We have thus obtained that
\beq
S\subset Q~~\mb{\rm and}~~S~\mb{bounded}\qu\Rightarrow\qu\LL^\m(S\cap\Tm_k)=0~.\lab{sq} \eeq
Using that every open subset of $\RM$ can be written as a countable union\footnote{From \eq{sq},
adjacent cubes can also intersect in our case.}
 of cubes in $\RM$ \cite{Rudin76,Evans92,Stein05}, the result in \eq{sq} entails that $\LL^\m(\Tm_\ka\cap Q)=0$.

Next, $Q^c\!=\ker\lp(\Pvt-\Pvs\rp)$. By  H\ref{aa}, $\dim\ker\lp(\Pvt\!-\Pvs\rp)\leq\m-1$.
Hence $\LL^\m(\Tm_\ka\cap Q^c)=0$.

The proof follows from the equality
$\LL^\m(\Tm_\ka)=\LL^\m(\Tm_\ka\cap Q)+\LL^\m(\Tm_\ka\cap Q^c)$.
\endproof

We exhibit a {\em closed negligible} subset of data in $\RM$ that can still meet
the equality in \eq{hq}.

\BP\lab{cor}
For $\be>0$ and $\k\in\Im$, put
\beq\Si_\k\deq\bigcup_{n=-\k}^{\k}~\bigcup_{\vt\in\Ombk}~
\bigcup_{\vs\in\Ombk}\lp\{g\in\RM~:~
\vt\neq\vs~~\mb{\rm and}~~g^T\Big(\Pvt-\Pvs\Big)g=n\be\rp\}~,\lab{S}
\eeq
where $\Ombk$ is given in~\eq{ga}.
Let H\ref{aa} hold for $\k$.
Then $\Si_\k$ is closed in $\RM$ and $\LL^\m\lp(\Si_\k\rp)=0$.
\EP

\proof
For some $n\in\{-\k,\cdots,\k\}$ and $(\vt,\vs)\in (\Ombk\x\Ombk)$ such that $\vt\neq\vs$, put
\[\Si\deq\lp\{g\in\RM:g^T\Big(\Pvt-\Pvs\Big)g=n\be\rp\}~.
\]
If $\#\vt\neq\#\vs$, then $\rank\big(\Pvt-\Pvs\big)\geq1$.
If $\#\vt=\#\vs$,  H\ref{aa} guarantees that  $\rank\big(\Pvt-\Pvs\big)\geq1$,
yet again. The number $n\be\in\RR$ is given.
According to \lem{H}, $\Si$ is a closed subset of $\RM$ and $\LL^\m\lp(\Si\rp)=0$.
The conclusion follows from the fact that $\Si_\k$ is a finite union of subsets like~$\Si$.
\endproof

{\em We assume hereafter that if H\ref{aa} holds for some $\k\in\Im$, data $d$ satisfy
\[d\in \{g\in\RM~:~g\not\in\Si_\k\}=\RM\sm\Si_\k~.\]}

\subsection{The unique global minimizer of $\bm{\Fd}$ is $\bm\k$-sparse for $\bm{\k\leq(\m-1)}$}\lab{kr}

We are looking for guarantees that $\Fd$ has a {\em unique} global minimizer $\hu$ obeying
\[ \|\hu\|_0\leq\k~~\mb{for some fixed}~~ \k\in\Im~.\]
This is the aim of the next theorem.
\BT \lab{Uglob}  Given $\k\in\Im$, let H\ref{aa} hold for $\k$,
$\be>\be_\k$ where $\be_\k$ meets \prop{ie} and $\Si_\k\subset\RM$
reads as in \eq{S}.
Consider that \[d\in\RM\sm\Si_\k~.\]
Then
\bit\item[\rm(i)] the set $\RM\sm\Si_\k$ is open and dense in $\RM$;
\item[\rm(ii)] $\F_d$ has a {\em unique global minimizer} $\hu$, and $\|\hu\|_0\leq\k$.
\eit
\ET

\proof Statement (i) follows from \prop{cor}.

Since $\be>\be_\k$,  all global
minimizers of $\Fd$ have their support in $\Ombk$ (\prop{ie}).
Using the fact that $d\in\RM\sm\Si_\k~$, the definition of $\Si_\k$ in \eq{S} shows that
\beq-\k\leq n\leq\k~~\mb{\rm and}~~(\vt,\vs)\in (\Ombk\x\Ombk),~\vt\neq\vs~~~\Rightarrow~~~
d^T\Big(\Pvt-\Pvs\Big)d\neq n\be~.\lab{gc}\eeq
The proof is conducted by contradiction.
Let $\hu$ and $\bu\neq\hu$ be two global minimizers of $\Fd$.
Then
\[\hsi\deq\si(\hu)\in\Ombk~~~\mb{\rm and}~~~\bsi\deq\si(\bu)\in\Ombk~,\]
and $\hsi~\neq~\bsi$.
By $\F_d(\hu)=\F_d(\bu)$,~\eq{hq} yields
\beq d^T\lp(\Phs-\Pbs\rp)d=\beta(\#\hsi-\#\bsi)~.\lab{qw}\eeq
An enumeration of all  possible values of $\#\hsi-\#\bsi$ shows that
\[\be\lp(\#\hsi-\#\bsi\rp)=n\be~~~\mb{for some}~~~n\in\{-\k,\cdots,\k\}~.\]
Inserting this equation into~\eq{qw} leads to
\[d^T\lp(\Phs-\Pbs\rp)d=n\be~~~\mb{for some}~~~n\in\{-\k,\cdots,\k\}~.\]
The last result contradicts~\eq{gc}; hence it violates the assumptions H\ref{aa} and $d\in\RM\sm\Si_\k$.
Consequently, $\Fd$ cannot have two global minimizers.
Since $\Fd$  always has global minimizers (\thm{gs}(i)), it follows that
$\Fd$ has a unique global minimizer, say $\hat u$.
And $\|\hat u\|_0\leq \mathsf{K}$ because $\sigma(\hat u)\in\overline{\Omega}_\mathsf{K}$.
\endproof

{\em For $\be>\be_\k$, the objective $\Fd$ in~\eq{fd} has a unique global minimizer and it is $\k$-sparse for
$\k\leq\m-1$.
For all $\;\k\in\Im$, the claim holds true in a generic sense.}
This is the message of \thm{Uglob} using \defi{ps}.


\section{Numerical illustrations}\lab{nm}
\subsection{On assumption H\ref{aa}}
Assumption H\ref{aa} requires that $\Pvt\neq\Pvs$ when $(\vt,\vs)\in\Om_r\x\Om_r$,
$\vt\neq\vs$ for all $r\leq\k\in\Im$.
From a practical viewpoint, the magnitude of $\lp(\Pvt-\Pvs\rp)$ should be discernible.
One way to assess the viability of H\ref{aa} for a matrix $A$ and $\k\in\Im$ is to calculate
\beqn\xi_\k(A)&\deq&\min_{r\in\Ik}~~\mu_r(A)~,\lab{xi}\\
\mb{where}~~~~~
 \mu_r(A)&=&\disp{\min_{\small\barr{c}(\vt,\vs)\in\Om_r\x\Om_r\\\vt\neq\vs\earr}
\lp\|\Pvt-\Pvs\rp\|_2,~~~\all r\in\Ik~.}\nn\eeqn
In fact, $\lp\|\Pvt-\Pvs\rp\|_2=\sin(\th)$, where $\th\in[0,\pi/2]$ is the maximum angle
between $\ran\lp(\Avt\rp)$ and $\ran\lp(\Avs\rp)$; see \cite[p. 456]{Meyer00}.
These subspaces have the same dimension and $\Pvt\neq\Pvs$ when
$(\vt,\vs)\in\Om_r\x\Om_r$, $\vt\neq\vs$ and $r\in\Ik$, hence $\th\in(0,\pi/2]$.
Consequently,
\[\mb{H\ref{aa} holds for}~~\k\in\Im~~~\Rightarrow~~~\mu_r(A)\in(0,1]~~\all r\in\Ik~~~\Rightarrow~~~\xi_\k(A)\in(0,1]~.\]
According to \eq{xi}, we have $\xi_\k\geq\xi_{\k+1}$, $\all \k\in\II_{\m-2}$.
Our guess that {\em assumption~H\ref{aa} is lightened when
$\k$ decreases} (see the comments following the proof of \thm{usg}) means that
\beq\xi_1(A)>\cdots>\xi_{\m-1}(A)~.\lab{tya}\eeq

We provide numerical tests on two subsets of real-valued random matrices for $\m=5$ and $\n=10$,
denoted by $\A^N_{20}$ and $\A^U_{1000}$.
The values of $\xi_\k(\cdot)$, $\k\in\Im=\II_4$, for every matrix in $\A^N_{20}$ and in
$\A^U_{1000}$, were calculated using an {\em exhaustive combinatorial search}.
{\em All tested matrices satisfy assumption H\ref{aa}}, which confirms \thm{usg} and its
consequences.
In order to evaluate the extent of H\ref{aa}, we computed the
{\em worst} and the {\em best} values of $\xi_\k(\cdot)$ over these sets:
\beq\left\{\barr{rcl}\xi_\k^{\mb{\FNS\it worst}}&=&\disp{\min_{A\in\A}\xi_\k(A)}~\\~\\
\xi_\k^{\mb{\FNS\it best}}&=&\disp{\max_{A\in\A}\xi_\k(A)}~\earr\right.~~~~~~\all\k\in\Im~,~~~
\A\in\lp\{\A^N_{20},~\A^U_{1000}\rp\}~.\lab{be}\eeq

\paragraph{Set $\A^N_{20}$}
This set was formed from 20 matrices $A^n$, ~$n\in\II_{20}$ of size $5\x 10$.
The components of each matrix $A^n$ were independent and uniformly drawn
    from the standard normal distribution with mean zero and variance one.
The values of $\xi_\k(\cdot)$ are depicted in Fig.~\ref{TH}.
We have\footnote{
This is why on the figure, in columns 10 and 17,
the green ``$\circ$'' and the red ``$\lozenge$'' overlap.}
$\xi_1(A^{10})=\xi_2(A^{10})$ and $\xi_1(A^{17})=\xi_2(A^{17})$.
In all other cases \eq{tya} is satisfied.
Fig.~\ref{TH} clearly shows that  $\xi_\k(\cdot)$ increases
as $\k$ decreases (from $\m-1$ to $1$).

\setlength{\unitlength}{1mm}
\bfig[h!] \bc
\bpic(143,81)
\put(0,1){\epsfig{figure=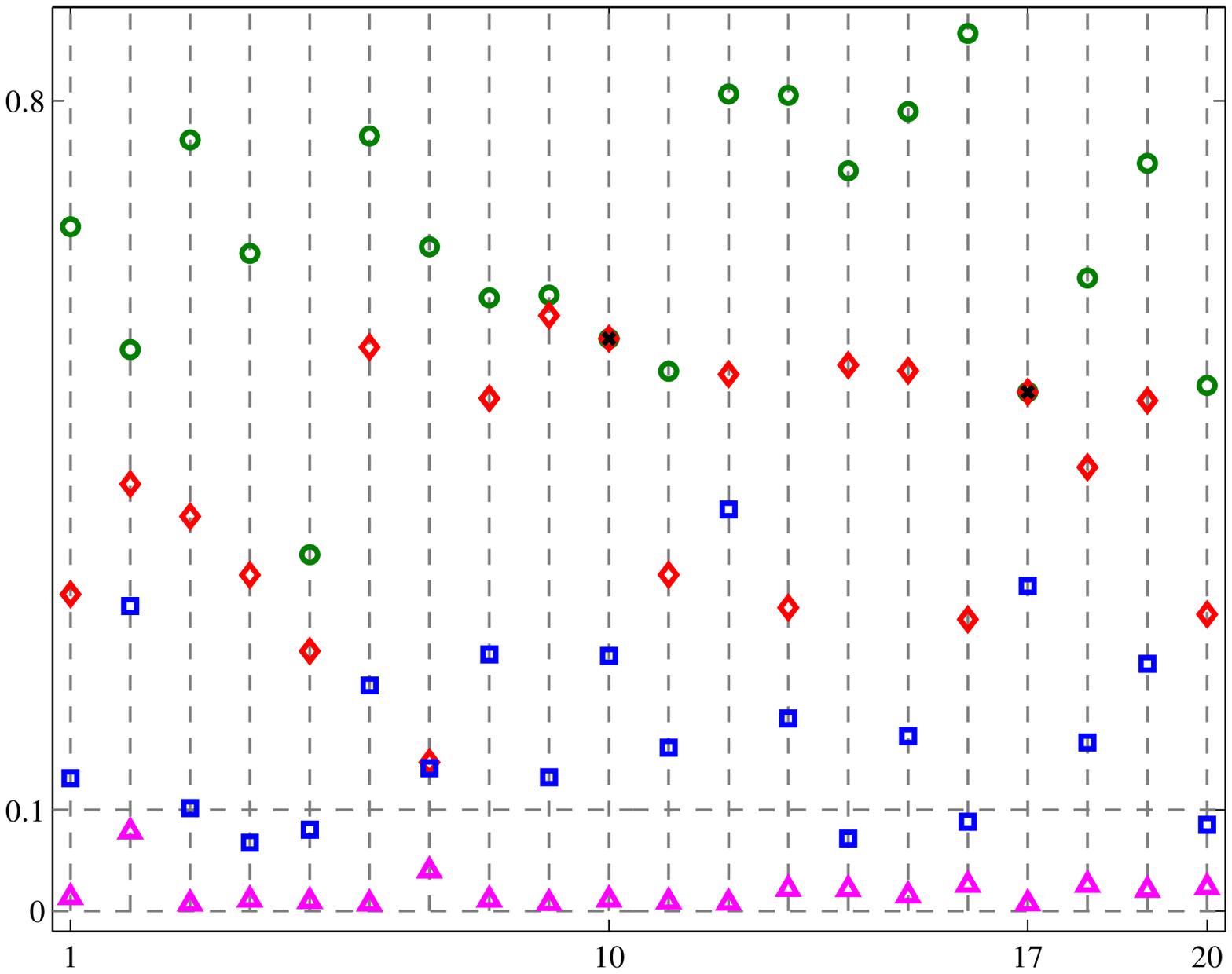,width=10.0cm,height=8cm}}
\put(20,-1.4){(a)~ $\xi_\k(\!A^n\!)$ for the matrices in $\A_{20}$}
\put(89,0.5){$n$}
\put(-5.6,61){\small$\xi_\k(\!A^n\!)$}
\put(101,1.2){\epsfig{figure=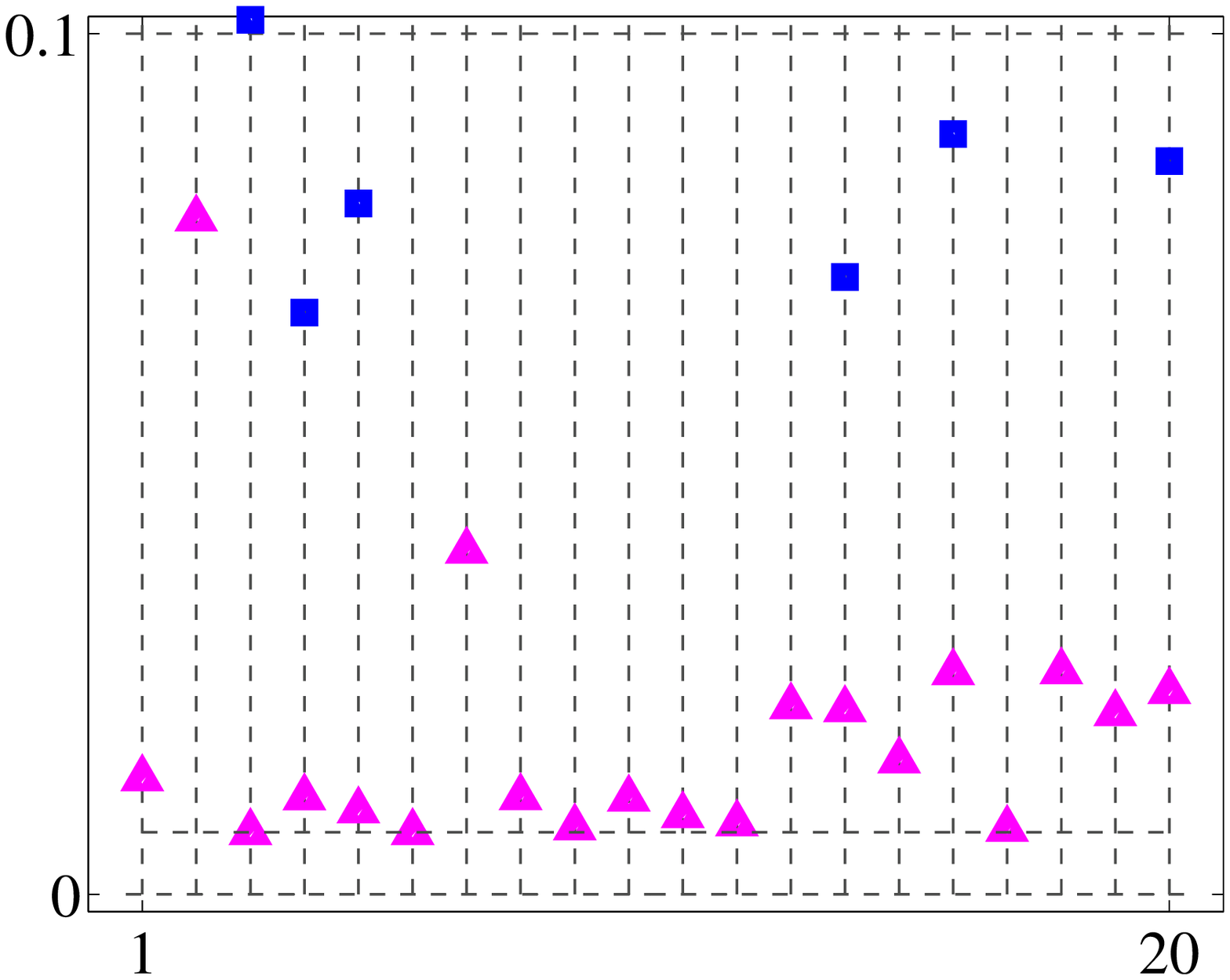,width=4.7cm,height=5cm}}
\put(101.5,30){z}
\put(101.5,28){o}
\put(101.5,26){o}
\put(100.8,24){m}
\put(107,-1.3){\small (b) Zoom of (a) -- $y$-axis}
\put(107,70){$\k=1$: \green{\bf green $\bm{\circ}$}}
\put(107,65){$\k=2$: \red{\bf red $\bm{\lozenge}$}}
\put(107,60){$\k=3$: \blue{\bf blue $\bm{\square}$}}
\put(107,55){$\k=4$: \mag{\bf magenta $\bm{\vartriangle}$}}
\epic \ec
\caption{$x$-axis: the list of the $20$ random matrices in $\A^N_{20}$.
(a) $y$-axis: the value $\xi_\k(A^n)$ according to~\eq{xi} for all $\k\in\Im$ and for all $n\in\II_{20}$.
The plot in (b) is a zoom of~(a) along the $y$-axis.}
\label{TH}
\efig

The worst and the best values of $\xi_\k(\cdot)$ over the whole set
$\A^N_{20}$ are displayed in Table~\ref{mmu}.
\btabe[h!]
\caption{The worst and the best values of $\xi_\k(A)$, for $\k\in\Im$, over the set $\A^N_{20}$,
see \eq{be}.}
\centering\mb{\begin{tabular}{l||c|c|c|c|}
\hline
& $\k=1$    &   $\k=2$  &   $\k=3$  &   $\k=4$\\
\hline
 $\disp{\xi_\k^{\mb{\FNS\it worst}}}$ & 0.3519   & 0.1467 &   0.0676  &  0.0072\\
\hline
 $\disp{\xi_\k^{\mb{\FNS\it best}}}$ &0.8666  &  0.5881  &  0.3966  &  0.0785\\

\hline
\end{tabular}}\lab{mmu}
\etabe

\paragraph{Set $\A^U_{1000}$}
The set $\A^U_{1000}$ was composed of one thousand $5\x 10$ matrices $A^n$, ~$n\in\II_{1000}$.
The entries of each matrix $A^n$ were independent and uniformly sampled on $[-1,1]$.
The obtained values for $\xi_\k^{\mb{\FNS\it worst}}$ and $\xi_\k^{\mb{\FNS\it best}}$,
calculated according to~\eq{be}, are shown in  Table~\ref{xmu}.

\btabe[h!]
\caption{The worst and the best values of $\xi_\k(A)$, for $\k\in\Im$, over the set $\A^U_{1000}$.}
\centering\mb{\begin{tabular}{l||c|c|c|c|}
\hline
& $\k=1$    &   $\k=2$  &   $\k=3$  &   $\k=4$\\
\hline
 $\disp{\xi_\k^{\mb{\FNS\it worst}}}$ & 0.1085 & 0.0235 &  0.0045 & 0.0001\\
\hline
 $\disp{\xi_\k^{\mb{\FNS\it best}}}$ &0.9526  &  0.8625  &  0.5379  &  0.1152\\
\hline
\end{tabular}}\lab{xmu}
\etabe
For $\k\in\II_3$, the {\em best} values of $\xi_\k(\cdot)$ were obtained for the same matrix, $A^{964}$.
Note that $\xi_4(A^{964})= 0.0425\gg\disp{\xi_4^{\mb{\FNS\it worst}}}$.
The {\em worst} values in Table~\ref{xmu} are smaller than those in
Table~\ref{mmu}, while the {\em best} values in Table~\ref{xmu} are larger than those in
Table~\ref{mmu}; one credible reason is that $\A^U_{1000}$ is much larger than $\A^N_{20}$.

\btabe[h!]
\caption{Percentage of the cases in $\A^U_{1000}$ when \eq{tya} fails to hold.}
\centering\mb{\begin{tabular}{c||c|c|c|}
\hline
& $\xi_1(A^n)=\xi_2(A^n)$&$\xi_2(A^n)=\xi_3(A^n)$&$\xi_3(A^n)=\xi_4(A^n)$\\
\hline
occurrences $\{n\}$ ~ & 5 $\%$&   1.6  $\%$&  0.1 $\%$\\
\hline
\end{tabular}} \lab{kx}
\etabe

Overall, \eq{tya} is satisfied on $\A^U_{1000}$---the
percentages in Table~\ref{kx} are pretty small.
All three tables and Figure~\ref{TH} agree with
our guess that~H\ref{aa} is more viable for smaller values of~$\k$.

{\em Based on the magnitudes for $\xi_\k^{\mb{\FNS\it best}}$ in Tables~\ref{mmu} and~\ref{xmu},
one can expect that there are some classes of matrices  (random or not) that
fit~H\ref{aa} for larger values of $\xi_\k(\cdot)$.}

\subsection{On the global minimizers of $\bm{\Fd}$}\lab{ong}

Here we summarize the outcome of a series of experiments corresponding to several
matrices $A\in\RMN$ where $\m=5$ and $\n=10$,
satisfying H\ref{aa} for $\k=\m-1$, different original vectors ${\ddot{u}}\in\RN$ and
data samples $d=A{\ddot{u}}+\mathrm{noise}$, for various values of $\be>0$.
In each experiment, we computed the complete list of all different strict (local) minimizers of $\Fd$,
say $\lp(\hu^i\rp)_{i=1}^n$.
Then the sequence of values $\lp(\Fd(\hu^i)\rp)_{i=1}^n$ was sorted in increasing order,
$\Fd\lp(\hu^{i_1}\rp)\leq\Fd\lp(\hu^{i_2}\rp)\leq\cdots\leq\Fd\lp(\hu^{i_n}\rp)~.$
A global minimizer $\hu^{i_1}$ is unique provided that
$\Fd\lp(\hu^{i_1}\rp)<\Fd\lp(\hu^{i_2}\rp)$.
In order to discard numerical errors, we also checked
whether $\lp|\Fd\lp(\hu^{i_1}\rp)-\Fd\lp(\hu^{i_2}\rp)\rp|$
is easy to detect.

{\em In all experiments we carried out, the following facts were observed:
\bit
\item The global minimizer
of $\Fd$ was unique---manifestly data $d$  never did belong to the {\em closed negligible} subset $\Si_\k$
in~\prop{cor}.
This confirms \thm{Uglob}.
\item
The global minimizers of $\Fd$ remained unchanged under large variations of~$\be$.
\item The necessary condition for a global minimizer in \prop{ou} was met.
\eit}

Next we present in detail two of these experiments where $\Fd$ is defined using
\begin{equation} A=\left[\begin{array}{c}
7~~~2~~~4~~~9~~~0~~~3~~~3~~~6~~~6~~~7\\
3~~~4~~~9~~~3~~~3~~~9~~~1~~~3~~~1~~~5\\
5~~~4~~~2~~~4~~~0~~~7~~~1~~~9~~~2~~~9\\
8~~~4~~~0~~~9~~~6~~~0~~~4~~~2~~~3~~~7\\
6~~~3~~~6~~~5~~~0~~~9~~~0~~~0~~~3~~~8
\end{array}\right]
~~\begin{array}{l}d=A{\ddot{u}}+n~,\\~\\
\mbox{~~~where $n$ is noise and}\\~\\
{\ddot{u}}=\big(\,0\,,~1\,,~8\,,~0\,,~3\,,~0\,,~0\,,~0\,,~0\,,~9\,\big)^T.\end{array}
\label{num}\end{equation}
Only integers appear in~\eq{num} for better readability.
We have $\rank(A)=\m=5$.
An exhaustive combinatorial test shows that
the arbitrary matrix $A$ in~\eq{num} satisfies~H\ref{aa} for $\k=\m-1$.
The values of $\xi_\k(A)$ are seen in Table~\ref{D}.
One notes that $\mu_2(A)>\mu_1(A)$; hence $\xi_1(A)=\xi_2(A)$.

\btabe[h!] \caption{The values of $\xi_\k(A)$ and $\mu_\k(A)$, $\all\k\in\Im$, for the matrix $A$
in~\eq{num}.}
\centering\mb{\begin{tabular}{l||c|c|c|c|}
\hline
& $\k=1$    &   $\k=2$  &   $\k=3$  &   $\k=4$\\
\hline
 $\xi_\k(A)$ & 0.2737 & 0.2737 &  0.2008 & 0.0564\\
\hline
\hline
 $\mu_\k(A)$ & 0.2737 & 0.2799 & 0.2008 & 0.0564\\
\hline
\end{tabular}} \lab{D}
\etabe

One {\em expects} (at least when data are noise-free) that the global minimizer $\hat u$ of $\mathcal{F}_d$ obeys
$\hat\sigma\subseteq\sigma({\ddot{u}})$, where ${\ddot{u}}$ is the original in~\eqref{num}, and that the
vanished entries of $\hat u$ correspond to the least entries of ${\ddot{u}}$.
This inclusion provides a partial way to rate the quality of the solution provided by a global minimizer
$\hat u$ of $\mathcal{F}_d$.

The experiments described hereafter correspond to two data samples relevant to
\eq{num}---without and with noise---and to several values of $\be>0$.

\paragraph{Noise-free data} The noise-free data in \eq{num} read as:
\beq d=A{\ddot{u}}=\big(~97\,,~~130\,,~~101\,,~~85\,,~~123~\big)^T.\lab{kli}\eeq
For different values of $\be$, the global minimizer $\hu$ is given in
Table~\ref{bem}.
\begin{table}[h!]
\caption{The global minimizer $\hat u$ of $\mathcal{F}_d$ and its value  $\mathcal{F}_d(\hat u)$
for the noise-free data $d$ in~\eqref{kli} for different values of $\beta$.
Last row: the original ${\ddot{u}}$ in \eqref{num}.}
\begin{tabular}{l}
\hline
\begin{tabular}{c|c|c|c|c}
$\beta$& The global minimizer $\hat u$ of $\mathcal{F}_d$ {\it ~~(row vector)}&$\|\hat u\|_0$ & $\mathcal{F}_d(\hat u)$ \\
\hline
$\begin{array}{c} 1\\ 10^2\\ 10^3\\ 10^4\\ 7\!\cdot\! 10^4
\end{array}$
&
$\begin{array}{cccccccccc}
0  & \mathbf{ 1}  &  \mathbf{8}     &  0  &  \mathbf{3}     &  0  &  0  &  0  &  0  &  \mathbf{9}\\
 0  &  0  & \mathbf{ 8.12}  &  0  &  \mathbf{3.31}  &  0  &  0  &  0  &  0  &  \mathbf{9.33}\\
 0 & 0 & 0 & 0 & 0 & \mathbf{12.58} & \mathbf{20.28} & 0 & 0 & 0\\
0 & \mathbf{29.95} & 0 & 0 & 0 & 0 & 0 & 0 & 0 & 0\\
0&0&0&0&0&0&0&0&0&0
\end{array}$&
$\begin{array}{c}4\\3\\2\\1\\0\end{array}$
&
$\begin{array}{c}
4\\ 301.52\\ 2179.3\\ 14144\\ 58864
\end{array}$
\\
\hline
\hline
\end{tabular}\\
$~~~~~~{\ddot{u}}~~~=~\,0~~~~~\;\mathbf{1}~~~~~~~~ \mathbf{8}~~~~~ 0~~~~\,\mathbf{3}~~~~~~~~ 0~~~~~~~~~ 0~~~~~~ 0~~~ 0~~~~~ \mathbf{9}$\\

\hline
\end{tabular}  \label{bem}
\end{table}
Since $\si({\ddot{u}})\in\Om$ and the data are noise-free,
$\Fd$ does not have global minimizers  with $\|\hu\|_0=5$.
Actually, applying \prop{ie} for $\tu={\ddot{u}}$ yields
$\be_{\m-1}=0$, hence for any $\be>0$ all global minimizers of $\Fd$ have a
support in $\Om=\Omb_{\m-1}$ (see \defi{om} and~\eq{ga}).
The global minimizer $\hu$ for $\be=1$ meets $\hu={\ddot{u}}$.
For $\be=100$, the global minimizer $\hu$ obeys $\hsi\deq\si(\hu)=\{3,5,10\}\subsetneqq\si({\ddot{u}})$ and
$\|\hu\|_0=3$---the least nonzero entry of the original ${\ddot{u}}$ is canceled, which is reasonable.
The global minimizers corresponding to $\be\geq 300$ are meaningless.
We could not find any positive value of $\be$ giving better $2$-sparse global minimizers.
Recalling that data are noise-free, this confirms \rem{x2}:
the global minimizers of $\Fd$ realize a only {\em pseudo}-hard thresholding.
For $\be\geq 7\cdot 10^4>\|d\|^2$, the global minimizer of $\Fd$ is $\hu=0$ which confirms~\rem{nu}.

\paragraph{Noisy data}
Now we consider {\em noisy data} in \eq{num} for
\beq  n=\big(~4\,,~~ -1\,,~~ 2\,,~~ -3\,,~~ 5~\big)^T~.\lab{noy}\eeq
This arbitrary noise yields a signal-to-noise ratio\footnote{
Let us denote $\ddot{d}=A{\ddot{u}}$ and $d=\ddot{d}+n$. The SNR reads \cite{Vetterli95}
$\mathrm{SNR}(\ddot{d},d)=10\,\log_{10}\frac{\sum_{i=1}^\m\lp(\ddot{d}[i]-\frac{1}{\m}\sum_{i=1}^\m \ddot{d}[i]\rp)^2}
{\sum_{i=1}^\m\lp(d[i]-\ddot{d}[i]\rp)^2}~.$
}
(SNR) equal to $14.07$ dB.
If $\be\leq0.04$, $\Fd$ has $252$ different strict global minimizers $\hu$
obeying $\|\hu\|_0=\m$ and $\Fd(\hu)=\be\m$ (recall \prop{bad}).
For $\be\geq0.05$, the global minimizer $\hu$ of $\Fd$ is
unique and satisfies $\si(\hu)\in\Om$.
It is given in Table~\ref{ben} for several values of $\be\geq0.05$.
\begin{table}[h!]
\caption{The global minimizer $\hat u$ of $\mathcal{F}_d$ and its value $\mathcal{F}_d(\hat u)$
for noisy data given by \eqref{num} and \eqref{noy},
for different values of $\beta$.
Last row: the original ${\ddot{u}}$.}
\begin{tabular}{l}
\hline
\begin{tabular}{c|c|c|c}
$\beta$& The global minimizer $\hat u$ of $\mathcal{F}_d$ {\it ~~(row vector)}&$\|\hat u\|_0$ & $\mathcal{F}_d(\hat u)$\\
\hline
$\begin{array}{c} 1\\ 10^2\\ 10^3\\ 10^4\\ 7\!\cdot\!10^4
\end{array}$
&
$\begin{array}{cccccccccc}
0 & \mathbf{6.02} & \mathbf{2.66} & \mathbf{6.43} & 0 & \mathbf{6.85} & 0 & 0 & 0 & 0\\
0 & 0 & \mathbf{8.23} & 0 & \mathbf{2.3} & 0 & 0 & 0 & 0 & \mathbf{9.71}\\
 0 & 0 & \mathbf{8.14} & 0 & 0 & 0 & 0 & 0 & 0 & \mathbf{10.25}\\
0 & 0 & 0 & 0 & 0 & 0 & 0 & 0 & 0 & \mathbf{14.47}\\
0&0&0&0&0&0&0&0&0&0
\end{array}$
&$\begin{array}{c}4\\3\\2\\1\\0\end{array}$
&
$\begin{array}{c}
4.0436\\301.94\\2174.8\\ 14473\\ 60559
\end{array}$
\\
\hline
\hline
\end{tabular}\\
$~~~~~~{\ddot{u}}~~~=\,~0~~~~~\mathbf{1}~~~~~~~ \mathbf{8}~~~~~~~\, 0~~~~~~ \mathbf{3}~~~~~~
 0~~~~~ 0~~~ 0~~\, 0~~~~~~ \mathbf{9}$\\
\hline
\end{tabular}  \label{ben}
\end{table}
For $\be=1$, the global minimizer is meaningless. We could not find any positive
value of $\be$ yielding a better global minimizer with a $4$-length support.
For the other values of $\be$, the
 global minimizer $\hu$ meets $\hsi\deq\si(\hu)\varsubsetneqq\si({\ddot{u}})$,
and its vanished entries correspond to the least entries in
the original ${\ddot{u}}$.
For $\be=100$, the global minimizer seems to furnish a good approximation to ${\ddot{u}}$.
Observe that the last entry of the global minimizer $\hu[10]$,
corresponding to the largest magnitude in ${\ddot{u}}$, freely increases when
$\be$ increases from $10^2$ to $10^4$.
We tested a tight sequence of intermediate values of $\be$ without finding better results.
Yet again, $\be\geq 7\cdot 10^4>\|d\|^2$ leads to a unique null global minimizer 
(see \rem{nu}).

\setlength{\unitlength}{1mm}
\bfig[h!] \bc
\bpic(160,83)
\put(5,7){\epsfig{figure=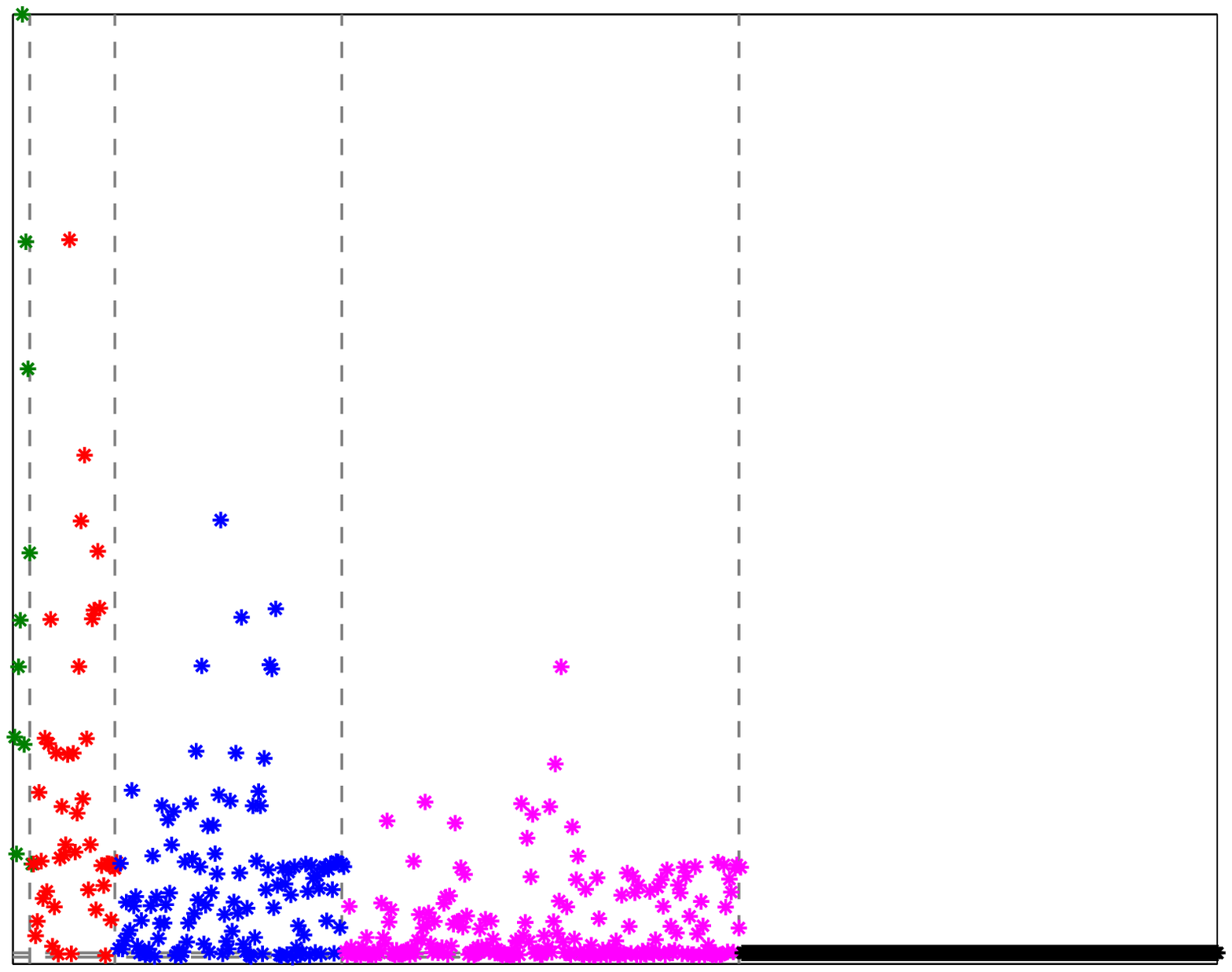,width=7.4cm,height=7.6cm}}
\put(85,8){\epsfig{figure=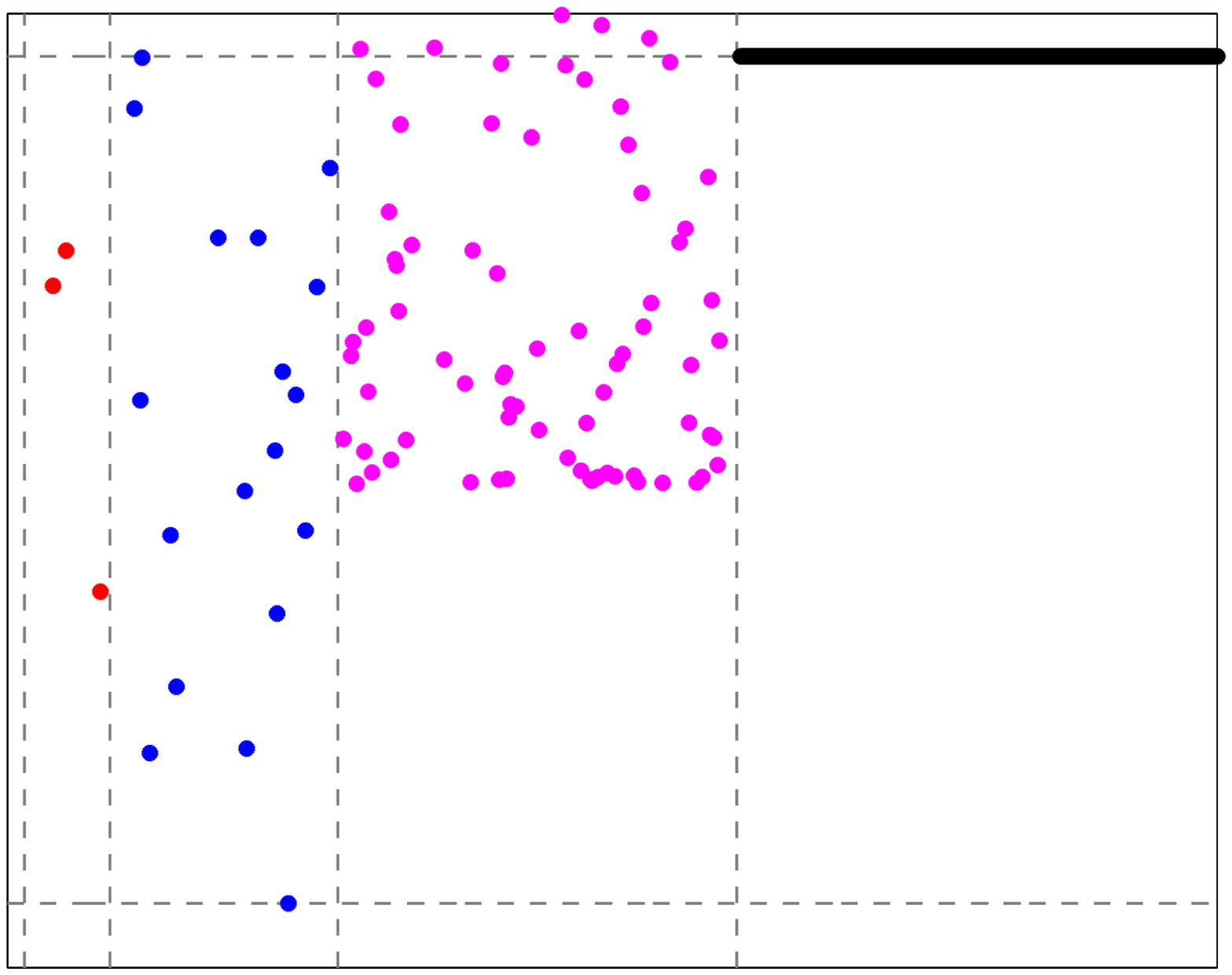,width=6.9cm,height=7.5cm}}
\put(6,0){\FNS (a) The values of all strict (local) minimizers of $\Fd$}
\put(94,1){\FNS (b) Zoom of (a) along the $y$-axis}
\put(-2,79){\FNS $4\;10^4$}
\put(0.5,4.5){\FNS $(0,0)$}
\put(75.3,4.5){\FNS $637$}
\put(65.3,4.5){\FNS $\{\hu\}$}
\put(0.0,7.5){\FNS $302$}
\put(-3,68){\FNS $\Fd(\hu)$}
\put(50,75){\FNS $\Fd(\hu)$ for}
\put(50,70){\FNS $\|\hu\|_0\!=\!1$: \green{\bf green}}
\put(50,65){\FNS $\|\hu\|_0\!=\!2$: \red{\bf red}}
\put(50,60){\FNS $\|\hu\|_0\!=\!3$: \blue{\bf blue}}
\put(50,55){\FNS $\|\hu\|_0\!=\!4$: \mag{\bf magenta}}
\put(50,50){\FNS $\|\hu\|_0\!=\!5$: {\bf black}}
\put(79.5,68){\FNS $\Fd(\hu)$}
\put(80.0,12){\FNS $302$}
\put(80.0,77.1){\FNS $500$}
\put(79.8,68){\FNS $\mathcal{F}_d(\hat u)$}
\put(149.2,5){\FNS $637$}
\put(140.5,5){\FNS $\{\hu\}$}
\put(87.5,5){\FNS $0$}
\put(91.9,18){\FNS $\|\hu\|_0\!=\!3$}
\put(109.7,18){\FNS $\|\hu\|_0=4$}
\put(133.3,18){\FNS $\|\hu\|_0=5$}
\epic \ec
\caption{All $\bm{638}$ strict (local) minima of $\Fd$ in \eq{num} for $\be=100$ and data $d$ corrupted
with the arbitrary noise in \eq{noy}.
The $x$-axis lists all strict (local) minimizers $\{\hu\}$ of $\Fd$ sorted according
to their $\ell_0$-norm $\|\hu\|_0$ in increasing order.
(a) The $y$-axis shows the value $\Fd(\hu)$ of these minimizers marked with a star.
The value of $\Fd$ for $\hu=0$ is not shown because it is too large ($\Fd(0)=60559=\|d\|^2$).
(b) A zoom of (a) along the $y$-axis.
It clearly shows that $\Fd$ has a very recognizable unique global minimizer}
\label{all}
\efig

Figure~\ref{all} shows the value $\Fd(\hu)$ of all the strict local minimizers of $\Fd$ for $\be=100$.
In the zoom in Figure~\ref{all}(b) it is easily seen that the global minimizer is unique (remember \thm{Uglob}).
It obeys $\|\hu\|_0=3$ and $\Fd(\hu)=301.94$.
One observes that $\Fd$ has $252=\#\Om_\m$ different strict local minimizers $\hu$
with $\|\hu\|_0=5=\m$ and $\Fd(\hu)=500=\be\m$.
This confirms \prop{bad}---obviously $d$ does not belong to the closed negligible
subset $\Qm_\m$ described in the proposition.

\section{Conclusions and perspectives}\lab{Con}

We provided a detailed analysis of the (local and global) minimizers of a
regularized objective $\Fd$ composed of a quadratic data fidelity term and an $\ell_0$ penalty
weighted by a parameter $\be>0$, as given in \eq{fd}.
We exhibited easy necessary and sufficient conditions ensuring that a (local) minimizer $\hu$ of $\Fd$ is strict
(\thm{ra}). The global minimizers of $\Fd$ (whose existence was proved) were shown to be strict as well (\thm{gs}).
Under very mild conditions, $\Fd$ was shown to have a unique global minimizer (\thm{Uglob}).
Other interesting results were listed in the abstract. Below we pose some
perspectives and open questions raised by this work.
\bit
\item The relationship between the value of the regularization parameter $\be$ and the
sparsity of the global minimizers of $\Fd$ (\prop{ie}) can be improved.

\item The {\em generic} linearity in data $d$ of each strict (local) minimizer of $\Fd$ (subsection~\ref{els})
should be exploited to better characterize the global minimizers of $\Fd$.

\item Is there a simple way to check whether assumption H\ref{aa} is satisfied
by a given matrix $A\in\RMN$ when $\mathsf{N}$ and  $\m<\n$ are large?
\rem{P} and in particular~\eq{yr} could help to discard some nonrandom matrices.
Conversely, one can ask whether there is a
systematic way to construct matrices $A$ that satisfy H\ref{aa}.

An alternative would be to exhibit families of matrices that satisfy H\ref{aa} for
 large values of $\xi_\k(\cdot)$, where the latter quantifiers
are defined in equation~\eq{xi}.

\item A proper adaptation of the results to matrices $A$ and data $d$ with complex entries should not present inherent
difficulties.

\item The theory developed here can be extended to MAP energies of the form evoked in~\eq{map}.
This is important for the imaging applications mentioned there.

\item Based on \cor{nor}, and Remarks~\ref{x2} and \ref{ifn},
and the numerical tests in subsection \ref{ong}, one
is justified in asking for conditions ensuring that the global minimizers of $\Fd$ perform
a valid work.
Given the high quality of the numerical results provided in many papers (see e.g.,~\cite{Robini07,Robini10}),
the question deserves attention.
\eit

There exist numerous algorithms aimed at approximating a (local) minimizer of $\Fd$.
As a by-product of our research, we obtained simple rules to verify whether or not an algorithm
could find
\bit\item[-] a (local) minimizer $\hu$ of $\Fd$---by checking whether $\hu$ satisfies \eq{ee} in \cor{nor};

\item[-]
and whether this local minimizer is strict by testing whether the submatrix whose
columns are indexed by the support of $\hu$ (i.e., $A_{\sigma(\hat u)}$) has full column rank (\thm{ra}).
\eit
Some properties of the minimizers of $\Fd$ given in this work can be inserted in numerical
schemes in order to quickly escape from shallow local minimizers.

Many existing numerical methods involve a studious choice of the regularization parameter  $\be$,
and some of them are proved to converge to a local minimizer of $\Fd$.
{\em We have seen that finding a (strict or nonstrict) local minimizer of $\Fd$
is easy and that it is independent of the value of $\be$ (Corollaries~\ref{nor} and \ref{blg}).
{\em It is therefore obscure what meaning to attach to
``choosing a good $\be$ and proving (local) convergence''.}}
Other successful algorithms are not guaranteed to converge to a local minimizer of $\Fd$.
Whenever algorithms do a good job, the choice of $\be$, the assumptions on $A$ and on $\|\hu\|_0$,
and the iterative scheme and its initialization {\em obviously} provide a tool for selecting a
meaningful solution by minimizing $\Fd$.
{\em There is a theoretical gap that needs clarification.}

{\em The connection between the existing algorithms and the description of the minimizers
exposed in this paper deserves deep exploration. What conditions ensure that an algorithm
minimizing $\Fd$ yields meaningful solutions? Clearly, showing local convergence does not answer
this important question.}

One can expect such research to give rise to innovative and more efficient algorithms
enabling one to compute relevant solutions by minimizing the tricky objective $\Fd$.

\section{Appendix}

\subsection{Proof of \lem{tl}}\lab{dtl}
 Since $\hu\neq0$, the definition of $\hsi$ shows that $\disp{\min_{i\in\hsi}\big|\,\hu[i]\,\big|>0}$.
Then  $\rho$ in~\eq{rho} fulfills $\rho>0$.

\paragraph{\rm(i)}~  Since $\#\hsi\geq1$, we have
\beqn i\in\hsi~,~~v\in\Bm_\infty(0,\rho)~~&\Rightarrow&~~\max_{j\in\hsi}\big|\,v[j]\,\big|<\rho\nn\\
&\Rightarrow&~~\max_{j\in\hsi}\big|\,v[j]\,\big|<\min_{j\in\hsi}\big|\,\hu[j]\,\big|\nn\\
&\Rightarrow&~~|\hu[i]+v[i]|\geq |\hu[i]|-|v[i]|\nn\\
&&~~\geq\min_{j\in\hsi} |\hu[j]|-\max_{j\in\hsi}|v[j]|\geq\rho-\max_{j\in\hsi}|v[j]|>0\nn\\
&\Rightarrow&~~
\hu[i]+v[i]\neq0~\nn\\
\Big[\mb{\rm by~\eq{phi}~}\Big]~~~~~~~~&\Rightarrow&~~\phi\lp(\hu[i]+v[i]\rp)=\phi\lp(\hu[i]\rp)=1~.
\lab{ry}\eeqn
If $\hsi^c=\void$ the result is proved. Let $\hsi^c\neq\void$.
Then $\hu[i]=0=\phi\lp(\hu[i]\rp)$, $\all i\in\hsi^c$.
Inserting this and~\eq{ry} into
\[\sum_{i\in\IN}\phi\big(\,\hu[i]+v[i]\,\big)=
\sum_{i\in\hsi}\phi\big(\,\hu[i]+v[i]\,\big)+\sum_{i\in\hsi^c}\phi\big(\,\hu[i]+v[i]\,\big)~\]
proves claim (i).

\paragraph{\rm(ii)}~
Using the fact that $\|A(\hu+v)-d\|^2=\|A\hu-d\|^2+\|Av\|^2+2\<Av,\,A\hu-d\>$, one obtains
\beqn v\in\Bm_\infty(0,\rho)\sm\Km_{\hsi}~~~\Rightarrow~~~
\Fd(\hu+v)&=&\|A\hu-d\|^2+\|Av\|^2+2\<Av,\,A\hu-d\>\nn\\
\Big[\mb{\rm by \lem{tl}(i)}\Big]~~~~&&+\be\sum_{i\in\hsi}\phi\lp(\hu[i]\rp)+\be\sum_{i\in\hsi^c}\phi\lp(v[i]\rp)\nn\\
\Big[\mb{\rm using~\eq{fds}~}\Big]~~~~&=&\Fd(\hu)+\|Av\|^2+2\<Av,\,A\hu-d\>+\be\sum_{i\in\hsi^c}\phi\lp(v[i]\rp)\nn\\
&\geq&\Fd(\hu)-\big|2\<v,\,A^T(A\hu-d)\>\big|+\be\|v_{\hsi^c}\|_0\nn\\
\Big[\mb{\rm by H\"{o}lder's inequality}\Big]~~~~&\geq&\Fd(\hu)-2\|v\|_\infty\,\|A^T(A\hu-d)\|_1+\be \|v_{\hsi^c}\|_0~.\lab{pop}
\eeqn
If $\#\hsi^c=0$, then $\Km_{\hsi}=\RN$, so $v\in\RR^0$ and $\|v\|_\infty=0$; hence we have the inequality.

Let  $\#\hsi^c\geq1$. For $v\not\in\Km_{\hsi}$, there at least one index $i\in\hsi^c$ such that $v[i]\neq0$;
hence  $\|v_{\hsi^c}\|_0\geq1$.
The definition of $\rho$ in~\eq{rho} shows that
\beqnn v\in\Bm_\infty(0,\rho)\sm\Km_{\hsi}~&\Rightarrow&~
    -\|v\|_\infty>-\rho\geq-\,\frac{\be}{2\Big(\|A^T(A\hu-d)\|_1+1\Big)}\\
~&\Rightarrow&~-2\|v\|_\infty\,\|A^T(A\hu-d)\|_1+\be\|v_{\hsi^c}\|_0
>-\,\frac{2\be\|A^T(A\hu-d)\|_1}{2\Big(\|A^T(A\hu-d)\|_1+1\Big)}+\be>0~.
\eeqnn
Introducing the last inequality into~\eq{pop} shows that for $\#\hsi^c\geq1$, the inequality in (ii) is strict.

\subsection{Proof of \prop{ou}}\lab{pou}
If $\hu=0$, the statement is obvious. We focus on $\hu\neq0$.
For an arbitrary $i\in\IN$, define
\[\hu^{(i)}\deq\big(\hu[1],\cdots,\hu[i-1],\,0,\,\hu[i+1],\cdots,\hu[\n]\big)\in\RN~.\]
We shall use the equivalent formulation of $\Fd$ given in~\eq{fds}.
Clearly\footnote{Using the definition of $\hu^{(i)}$, we have
$\hu^{(i)}=A_{(\IN\sm\{i\})}\hu_{(\IN\sm\{i\})}~,$
hence $A\hu^{(i)}$ is independent of $\hu[i]$.
},
\[\Fd(\hu)=\Fd\big(\hu^{(i)}+e_i\hu[i]\big)=\|A\hu^{(i)}+a_i\hu[i]-d\|^2+
   \be \sum_{j\in\IN}\phi\lp(\hu^{(i)}[j]\rp)+\phi\lp(\hu[i]\rp)~.\]
Consider $f:\RR\to\RR$ as given below
\beq f(t)\deq\Fd\lp(\hu^{(i)}+e_it\rp). \lab{fi}\eeq
Since $\hu$ is a global minimizer of $\Fd$, for any $i\in\IN$, we have
\beqnn f\lp(\hu[i]\rp)&=&\Fd\big(\hu^{(i)}+e_i\hu[i]\big)\\
    &\leq&\Fd\big(\hu^{(i)}+e_i\,t\big)=f(t)~~~\all t\in\RR~.\eeqnn
Equivalently, for any $i\in\IN\,$, $\,f\lp(\hu[i]\rp)\,$ is the global minimum of
$f(t)$ on $\RR$.
Below we will determine the global minimizer(s) $\h t=\hu[i]$ of $f$ as given in~\eq{fi}, i.e.,
\[\h t=\hu[i]=\arg\min_{t\in\RR}f(t)~.\]
In detail, the function $f$ reads as
\beqn f(t)&=&\|A\hu^{(i)}+a_it-d\|^2+
   \be \sum_{j\in\IN}\phi\lp(\hu^{(i)}[j]\rp)+\be\phi(t)\nn\\
&=&\|A\hu^{(i)}-d\|^2+\|a_i\|^2t^2+2t\<a_i,A\hu^{(i)}-d\>
+\be \sum_{j\in\IN}\phi\lp(\hu^{(i)}[j]\rp)+\be\phi(t)\nn\\
&=&\|a_i\|^2t^2+2t\<a_i,A\hu^{(i)}-d\>+\be\phi(t)+C\lab{fin0}~,
\eeqn
where
\[C=\|A\hu^{(i)}-d\|^2+\be \sum_{j\in\IN}\phi\lp(\hu^{(i)}[j]\rp)~.\]
Note that $C$ does not depend on $t$.
The function $f$ has two local minimizers denoted by
$\h t_0$ and $\h t_1$.
The first is
\beq\h t_0=0~~~\Rightarrow~~~f(\h t_0)=C~.\lab{min1}\eeq
The other one,   $\h t_1\neq0$, corresponds to  $\phi(t)=1$. From~\eq{fin0}, $\h t_1$ solves
\[2\|a_i\|^2t+2\<a_i,A\hu^{(i)}-d\>=0~.\]
Recalling that $a_i\neq0$, $\all i\in\IN$ (see \eq{ao}), it follows that
\beq\h t_1=-\frac{\<a_i,A\hu^{(i)}-d\>}{\|a_i\|^2}~~~\Rightarrow~~~
    f(\h t_1)=-\frac{\<a_i,A\hu^{(i)}-d\>^2}{\|a_i\|^2}+\be+C~.
\lab{min2}\eeq
Next we check whether $\h t_0$ or $\h t_1$ is a global minimizer of $f$.
From~\eq{min1} and~\eq{min2} we get
\[ f(\h t_0)-f(\h t_1)=\frac{\<a_i,A\hu^{(i)}-d\>^2}{\|a_i\|^2}-\be~.
\]
Furthermore,
\beqn f(\h t_0)<f(\h t_1)~~&\Rightarrow&~~\hu[i]=\h t_0=0~,\nn\\
f(\h t_1)<f(\h t_0)~~&\Rightarrow&~~\hu[i]=\h t_1=-\frac{\<a_i,A\hu^{(i)}-d\>}{\|a_i\|^2}~,\lab{min3}\\
f(\h t_0)=f(\h t_1)~~&\Rightarrow&~~\mb{$\h t_0$ and $\h t_1$ are global minimizers of $f$.}\nn
\eeqn
In particular, we have
\beqn f(\h t_1)\leq f(\h t_0)~~&\Leftrightarrow&~~\<a_i,A\hu^{(i)}-d\>^2\geq\be\|a_i\|^2\lab{Delta1}\\
\Big[\mb{\rm by~\eq{min3}}\Big]~~~~~~~~&\Rightarrow&~
    |\hu[i]|=\frac{\lp|\<a_i,A\hu^{(i)}-d\>\rp|}{\|a_i\|^2}\nn\\
\Big[\mb{\rm by~\eq{Delta1}}\Big]~~~~~~~~
&&~~~~~~~~\geq\frac{\sqrt{\be}\|a_i\|}{\|a_i\|^2}=\frac{\sqrt{\be}}{\|a_i\|}~.\nn
\eeqn
It is clear that the conclusion holds true for any $i\in\IN$.

\subsection{ Proof of Proposition~\ref{pal}}\lab{palp}
The asymptotic function $\lp(\Fd\rp)_\infty(v)$ of $\Fd$ can be calculated
according to\footnote{In the nonconvex case, the notion of asymptotic functions and the representation
formula were first given by J.P. Dedieu \cite{Dedieu77}.}
~\cite[Theorem 2.5.1]{Auslender03}
\[ \lp(\Fd\rp)_\infty(v)=\liminf_{\barr{c}v'\to v\\ t\to\infty\earr}\frac{\Fd(tv')}{t}~.\]
Then
\beqnn \lp(\Fd\rp)_\infty(v)&=&\liminf_{\barr{c}v'\to v\\ t\to\infty\earr}\frac{\|Av'-d\|^2+\be\|v'\|_0}{t}\nn\\
&=&\liminf_{\barr{c}v'\to v\\ t\to\infty\earr}\lp(t\|Av'\|^2-2\<d,Av'\>+\frac{\|d\|^2+\be\|v'\|_0}{t}\rp)\nn\\~\nn\\
&=&\left\{\barr{lll} 0 & \mb{if} & v\in \ker(A)~,\\+\infty & \mb{if} & v\not\in \ker(A)~.\earr
\right.\lab{auv}
\eeqnn
Hence
\beq\ker\lp((\Fd)_\infty\rp)=\ker(A)~,\lab{wd}\eeq
where $\ker\lp((\Fd)_\infty\rp)=\{v\in\RN~:~(\Fd)_\infty(v)=0\}$.

Let $\{v_k\}$ satisfy~\eq{aal} with $v_k\,\|v_k\|^{-1}\to \bv\in\ker(A)$.
Below we compare the numbers $\|v_k\|_0$ and $\|v_k-\rho \bv\|_0$ where $\rho>0$.
There are two options.
\ben
\item Consider that $i\in\si(\bv)$, that is, $\disp{\bv[i]=\lim_{k\to\infty} v_k[i] \, \|v_k\|^{-1}\neq 0}$.
Then $|\,v_k[i]\,|>0$ for all but finitely many $k$ as otherwise, $v_k[i] \, \|v_k\|^{-1}$
would converge to 0.
Therefore, there exists $k_i$ such that
\beq \lp|\,v_k[i]-\rho\,\bv[i]\,\rp|\geq0\qu\mb{and}\qu |\,v_k[i]\,|>0 \qu\all k\geq k_i~.\lab{wb}\eeq

\item If $i\in(\si(\bv))^c$, i.e. $\bv[i]=0$, then clearly
\beq v_k[i]-\rho\, \bv[i]=v_k[i] ~.\lab{wc}\eeq
\een
Combining \eq{wb} and \eq{wc}, the definition of $\|\cdot\|_0$ using $\phi$ in \eq{phi} shows that
\beq\|v_k-\rho\,\bv\|_0\leq\|v_k\|_0 \qu\all k\geq k_0\deq\max_{i\in\si(\bv)}k_i~. \lab{we}\eeq

By \eq{wd}, $A\bv=0$. This fact, jointly with~\eq{we}, entails that
\beqn\Fd(v_k-\rho\bv)&=&\|A(v_k+\rho\bv)-d\|^2+\be\|v_k-\rho\bv\|_0\nn\\
&=&\|Av_k-d\|^2+\be\|v_k-\rho\bv\|_0\nn\\
&\leq& \|Av_k-d\|^2+\be\|v_k\|_0=\Fd(v_k)\qu\all k\geq k_0~.\nn\lab{hb}\eeqn
It follows that for any $k\geq k_0$ we have
\[v_k\in\lev(\Fd,\la_k)\qu \Rightarrow\qu v_k-\rho\bv\in\lev(\Fd,\la_k)~,\]
and thus $\Fd$ satisfies \defi{als}.

\subsection{Proof of \prop{ie}}\lab{pie}
Given $\k\in\II_{\m-1}$, set
\beq\Um_{\k+1}\deq\bigcup_{\vt\subset\IN}\lp\{\bu~:~\bu~~\mb{solves \cpa}~~\mb{\rm and}~~
\|\bu\|_0\geq\k+1\rp\}~.\lab{ws}\eeq
\bit\item
Let $\Um_{\k+1}\neq\void$.
By \prop{cp}, for any $\be>0$, $\Fd$ has a (local) minimum at
each $\bu\in\Um_{\k+1}$. Thus
\beq\bu~~\mb{is a  (local) minimizer  of $\Fd$ and}~~\|\bu\|_0\geq\k+1
~~~\Leftrightarrow~~~\bu\in\Um_{\k+1}~.\lab{koi}\eeq
Then for any $\be>0$
\beq\Fd(\bu)\geq\be(\k+1)\qu\all\bu\in\Um_{\k+1}~.\lab{kf}\eeq

Let $\tu$ be defined by\footnote{Such a $\tu$ always exists; see subsection~\ref{ntd}.
By \prop{cp} and \thm{ra}, it is uniquely defined.}:
\[\tu~~\mb{solves~\cpa~for~some}~~\vt\in\Om_\k~.\]
Then
\beq \|\tu\|_0\leq\k~\lab{ku}~.\eeq
Set $\be$ and $\be_\k$ according to
\beq\be>\be_\k\deq\|A\tu-d\|^2~.\lab{kn}\eeq
For such a $\be$ we have
\beqnn\Fd(\tu)&=&\|A\tu-d\|^2+\be\|\tu\|_0\\
\Big[\mb{\rm by~\eq{ku}~and~\eq{kn}~}\Big]&<&\be+\be\k=\be(\k+1)\\
\Big[\mb{\rm by~\eq{kf}~}\Big]&\leq&\Fd(\bu)~~~~\all\bu\in\Um_{\k+1}~.
\eeqnn
Let $\hu$ be a global minimizer of $\Fd$. Then
\[\Fd(\hu)\leq\Fd(\tu)<\Fd(\bu)~~~~\all\bu\in\Um_{\k+1}~.\]
Using \eq{ws}-\eq{koi}, we find $$\|\hu\|_0\leq\k~.$$
\item
$\Um_{\k+1}=\void$ entails that\footnote{Let $A=(e_1,e_2,e_3,e_4,e_1)\in\RR^{4\x 5}$ and $d=e_1\in\RR^4$.
For $\k=\m-1=3$ one can check that $\Um_{\k+1}=\void$.}
\beq\bu~~\mb{solves~\cpa~for}~~\vt\subset\IN,~~\#\vt\geq\k+1~~~\Rightarrow~~~\|\bu\|_0\leq\k~.\lab{ha}\eeq
Let $\hu$ be a global minimizer of $\Fd$. By \eq{ha} we have
\[\|\hu\|_0\leq\k~.\]
\eit
According to \thm{gs}(ii), any global minimizer of $\Fd$ is strict, hence
$\disp{\si(\hu)\in\Ombk}~.$
\section*{Acknowledgments} The author would like to thank the anonymous Reviewer $2$
for very helpful remarks and suggestions.

\def\ieeeASSP{IEEE Transactions on Acoustics Speech and Signal Processing}
\def\ieeeIP{IEEE Transactions on Image Processing}
\def\ieeeIT{IEEE Transactions on Information Theory}
\def\ieeePAMI{IEEE Transactions on Pattern Analysis and Machine Intelligence}
\def\ieeeSP{IEEE Transactions on Signal Processing}
\def\JRSSB{Journal of the Royal Statistical Society B}
\def\CRAS{Compte-rendus de l'acad\'emie des sciences}
\def\BMK{Biometrika}
\def\IJCV{International Journal of Computer Vision}
\def\siamMMS{SIAM Journal on Multiscale Modeling and Simulation}
\def\JAS{Journal of Applied Statistics}
\def\siamR{SIAM Review}
\def\siamIS{SIAM Journal on Imaging Sciences}
\def\JGO{Journal of Global Optimization}
\def\JSC{Journal of Scientific Computing}
\def\MC{Mathematics of Computation}
\def\MP{Mathematical Programming}
\def\JFAA{Journal of Fourier Analysis and Applications}

\def\up{\uppercase}


\end{document}